\documentclass[toc=bib,parskip=half]{scrartcl}
\usepackage[utf8]{inputenc}
\usepackage{amsmath}
\usepackage{amsfonts}
\usepackage{amsthm}
\usepackage{amssymb}
\usepackage{thmtools}
\usepackage{enumerate}
\usepackage{hyperref}
\usepackage{enumitem}
\usepackage{stmaryrd}
\usepackage{mathrsfs}

\usepackage{graphicx}
\usepackage[dvipsnames]{xcolor}

\usepackage{tikz}
\usepackage{verbatim}
\usepackage{kbordermatrix}
\usepackage{todonotes}

\title{Schottky presentations of positive representations}
\author{Jean-Philippe Burelle \thanks{Burelle gratefully acknowledges support from the Natural Sciences and Engineering Research Council of Canada (NSERC), U.S. National Science Foundation grants DMS 1107452, 1107263, 1107367 ``RNMS: Geometric Structures and Representation Varieties'' (the GEAR Network). This project has received funding from the European Research Council (ERC) under the European Union’s Horizon 2020 research and innovation programme (ERC starting grant DiGGeS, grant agreement No 715982).} \and Nicolaus Treib \thanks{Treib gratefully acknowledges support from the Klaus Tschira Foundation, the European Research Council under ERC-Consolidator grant 614733 and ERC starting grant DiGGeS, grant agreement No 715982, the RTG 2229 grant of the German Research Foundation, and the GEAR Network.}}

\AtBeginDocument{%

}

\declaretheoremstyle[spaceabove=\topsep,spacebelow=0pt,bodyfont=\normalfont]{scdef}
\declaretheoremstyle[spaceabove=\topsep,spacebelow=0pt,bodyfont=\itshape]{scthm}
\declaretheoremstyle[spaceabove=\topsep,spacebelow=0pt,headfont=\normalfont\itshape,notefont=\normalfont\itshape,notebraces={}{},headformat={\NAME\NOTE},postheadspace=1em,qed=\qedsymbol]{scprf}

\declaretheorem[style=scthm,numberwithin=section,name=Theorem,    refname={Theorem,Theorems},        Refname={Theorem,Theorems}]        {Thm}
\declaretheorem[style=scthm,sharenumber=Thm, name=Lemma,      refname={Lemma,Lemmas},            Refname={Lemma,Lemmas}]            {Lem}
\declaretheorem[style=scthm,sharenumber=Thm, name=Corollary,  refname={Corollary,Corollaries},   Refname={Corollary,Corollaries}]   {Cor}
\declaretheorem[style=scthm,sharenumber=Thm, name=Proposition,refname={Proposition,Propositions},Refname={Proposition,Propositions}]{Prop}
\declaretheorem[style=scdef,sharenumber=Thm, name=Definition, refname={Definition,Definitions},  Refname={Definition,Definitions}]  {Def}
\declaretheorem[style=scdef,sharenumber=Thm, name=Remark,     refname={Remark,Remarks},          Refname={Remark,Remarks}]          {Rem}

\declaretheorem[style=scdef,sharenumber=Thm, name=Example,      refname={Example,Examples},            Refname={Example,Examples}]            {Ex}
\declaretheorem[style=scprf,unnumbered,name=Proof]{Prf}

\newtheoremstyle{TheoremNum}
    {\topsep}{\topsep}              
    {\itshape}                      
    {}                              
    {\bfseries}                     
    {.}                             
    { }                             
    {\thmname{#1}\thmnote{\bfseries #3}}
\theoremstyle{TheoremNum}

\numberwithin{equation}{subsection}

\DeclareMathOperator{\Sp}{Sp}

\DeclareMathOperator{\SO}{SO}

\DeclareMathOperator{\PSL}{PSL}
\DeclareMathOperator{\SL}{SL}
\DeclareMathOperator{\pr}{\mathbb{P}}
\DeclareMathOperator{\spr}{\mathbb{S}}

\newcommand{\RP}[1]{\mathbb{RP}^{#1}}
\newcommand{\bR}{\mathbb{R}}
\newcommand{\bN}{\mathbb{N}}

\newcommand{\bH}{\mathbb{H}}

\newcommand{\lift}[1]{\widehat{#1}}

\newcommand{\ival}[2]{(\!(#1,#2)\!)}
\newcommand{\sfival}{\ival{F_e}{F_{w_0}}}
\newcommand{\comp}{\mathcal{C}}

\newcommand{\dbfrac}[2]{\frac{\ #1\ }{\ #2\ }}
\newcommand{\chinv}[2]{\tau(#1,#2)}

\newcommand{\minor}[3]{{#1}^{#2}_{#3}}
\newcommand{\submat}[3]{#1\big[ \, \substack{#2\\#3} \, \big]}
\newcommand{\mind}[2]{\mathcal{I}(#1,#2)}

\newcommand{\Flag}[1]{\mathrm{Flag}(\bR^{#1})}
\newcommand{\oFlag}[1]{\mathrm{Flag}^+(\bR^{#1})}
\newcommand{\loFlag}[1]{\widehat{\mathrm{Flag}}^{\raisebox{-3.5pt}{$\scriptstyle +$}}(\bR^{#1})}
\newcommand{\equalo}{\overset{+}{=}}

\newcommand{\Gro}{\mathrm{Gr}^+}

\newcommand{\ltu}{\mathcal{L}}
\newcommand{\ltup}{\mathcal{L}^{>0}}
\newcommand{\ltun}{\mathcal{L}^{\geq 0}}
\newcommand{\utu}{\mathcal{U}}
\newcommand{\utup}{\mathcal{U}^{>0}}
\newcommand{\utun}{\mathcal{U}^{\geq 0}}

\newcommand{\bdry}{{\partial_\infty\Gamma}}

\newcommand{\Half}{\mathcal{H}}

\newcommand{\AdS}{\mathsf{AdS}}

\setlist[enumerate,1]{label={(\roman*)}}

\newcommand{\cycle}[3]{\overrightarrow{#1#2#3}}

\renewcommand{\kbldelim}{(}
\renewcommand{\kbrdelim}{)}

\DeclareMathOperator{\Aut}{Aut}

\DeclareMathOperator{\stab}{Stab}
\DeclareMathOperator{\Span}{span}
\DeclareMathOperator{\Hom}{Hom}

\DeclareMathOperator{\Var}{S}

\begin{document}

\maketitle
\begin{abstract}
    We show that the notion of $3$-hyperconvexity on oriented flag manifolds defines a partial cyclic order. Using the notion of interval given by this partial cyclic order, we construct Schottky groups and show that they correspond to images of positive representations in the sense of Fock and Goncharov. We construct polyhedral fundamental domains for the domain of discontinuity that these groups admit in the projective space or the sphere, depending on the dimension.
\end{abstract}
\tableofcontents

\section{Introduction}
Let $\Sigma$ be an oriented surface of negative Euler characteristic, $\Gamma=\pi_1(\Sigma)$ its fundamental group, and $G$ a simple real Lie group. Higher Teichm\"uller spaces are open subsets of the character variety $\Hom(\Gamma,G)/\!/G$ exhibiting properties similar to that of the classical Teichm\"uller space. For example, all representations in these spaces are faithful and discrete.

The first examples of such spaces were discovered in \cite{hitchin} and are now known as Hitchin components. For $G=\PSL(n,\bR)$, these are the connected components of the character variety containing representations which are the composition of a discrete and faithful representation into $\PSL(2,\bR)$ with the irreducible representation $\rho:\PSL(2,\bR)\rightarrow\PSL(n,\bR)$. The structure of these components was investigated in the foundational papers \cite{labourie} and \cite{FockGoncharov}. In the latter, an analog of the Hitchin component, the space of \emph{positive representations}, was defined for surfaces with boundary. Already in this early work, the importance of the cyclic structure on the boundary of the universal cover of a surface was clearly emphasized.

The second family of higher Teichm\"uller spaces to be investigated was that of \emph{maximal representations} in \cite{BIW}. Although the tools used to study maximal representations are generally different from those which were successful for Hitchin components and positive representations, the one feature which seems to connect these higher Teichm\"uller theories is the cyclic structure on the boundary of the group, and a compatible cyclic structure on a homogeneous space of $G$ (see \cite{BIW2} for details on this analogy). In \cite{BTSchottky}, the authors defined a notion of \emph{generalized Schottky group} of automorphisms of a space admitting a partial cyclic order, and showed that maximal representations are examples of this construction. This characterization was then used to build fundamental domains for the action of maximal representations into symplectic groups $\Sp(2n,\bR)$ on a domain of discontinuity in projective space.

The first goal of this paper is to define a partial cyclic order on the space of complete oriented flags in $\bR^n$, and show that positive representations into $\PSL(n,\bR)$ in the sense of Fock and Goncharov are generalized Schottky groups acting on this cyclically ordered space. The space of complete oriented flags is the quotient $\oFlag{n} := \PSL(n,\bR)/B_0$ of $\PSL(n,\bR)$ by the identity component of its Borel subgroup of upper triangular matrices. Its elements are sequences of nested subspaces in $\bR^n$ with a choice of orientation on each subspace, modulo the action of $-1$ if $n$ is even.

The partial cyclic order gives rise to a notion of intervals $I=\ival{F}{G}\subset\oFlag{n}$, which are open subsets of the oriented flag variety: The interval $I$ consists of all oriented flags $H$ such that $(F,H,G)$ is in cyclic configuration.
Each interval has an \emph{opposite}, obtained by reversing the endpoints, which we denote by $-I$. Generalized Schottky groups are then defined using a finite collection of intervals $I_j$ such that $I_j\subset -I_k$ whenever $j\neq k$. Each generator maps the opposite of some interval to another interval, analogously to the case of Schottky subgroups of $\PSL(2,\bR)$ acting on $\RP{1}$. 

\begin{Thm}
  Let $\Sigma$ be an oriented compact surface with boundary and let $\rho: \pi_1(\Sigma)\rightarrow \PSL(n,\bR)$ be a positive representation. Then, $\rho$ admits a presentation as a Schottky group pairing disjoint intervals in $\oFlag{n}$. Conversely, any Schottky group constructed this way using cyclically ordered intervals is the image of a positive representation.
\end{Thm}
\begin{Rem}
  In the converse part of the theorem, the dependence on the surface $\Sigma$ is hidden in the cyclic ordering of the intervals and the choice of pairings.
\end{Rem}
When the intervals chosen to define the Schottky group are not allowed to share endpoints, we call the resulting group \emph{purely hyperbolic}, again by analogy with the $\PSL(2,\bR)$ setting. In this case, we have a better understanding of the dynamics of the group. We show that the resulting representations are $B$-Anosov, where $B$ is a Borel subgroup. 

\begin{Thm}
  Let $F_g$ denote the free group on $g$ generators and let $\rho : F_g \rightarrow \PSL(n,\bR)$ be a purely hyperbolic Schottky representation. Then, $\rho$ is $B$-Anosov.
\end{Thm}
For $n=2k$ even, we then associate to an interval $I\subset \oFlag{2k}$ a \emph{halfspace} $\Half(I)\subset \RP{2k-1}$ bounded by a polyhedral hypersurface. This new notion of halfspace is related to previous constructions of fundamental domains for affine Schottky groups in dimension $3$ using surfaces called \emph{crooked planes}, introduced in \cite{DrummCrookedPlanes}. Crooked planes were generalized in various directions including to the setting of $3$-dimensional conformal Lorentzian geometry in \cite{FrancesCrookedSurf} and to $3$-dimensional anti-de Sitter geometry in \cite{DGKAntideSitter}. We explain how anti-de Sitter crooked halfspaces are related to $\RP{3}$ halfspaces in Appendix A. Interestingly, the groups for which these hypersurfaces bound fundamental domains in each case are not part of the class studied in this paper (positive representations). The ubiquity and effectiveness of (generalized) crooked planes in constructing fundamental domains for properly discontinuous actions of free groups remains mysterious in general.

Here are some key properties relating intervals and halfspaces:
\begin{itemize}
    \item Whenever $I\subset -J$, the halfspaces $\Half(I)$ and $\Half(J)$ are disjoint (\autoref{Thm:HalfspaceDisjointness});
    \item $\Half(-I) = \RP{2k-1}-\overline{\Half(I)}$, where the bar denotes closure in $\RP{2k-1}$ (\autoref{Lem:HalfSpaceComplementEven});
    \item If $I_1\supset I_2\supset\dots$ is a sequence of nested intervals with $\bigcap_j I_j=F\in \oFlag{2k}$, then $\bigcap \Half(I_j) = \mathbb{P}F^{(k)}$ (\autoref{Prop:HalfspacesAndIntervals}).
\end{itemize}
These properties allow us to build a fundamental domain in projective space $\RP{2k-1}$ by intersecting the complements of halfspaces associated to the disjoint intervals used in defining the Schottky group. In the case of purely hyperbolic (Anosov) representations, the orbit of this fundamental domain is the cocompact domain of discontinuity $D\subset \RP{2k-1}$ identified in \cite{GWDoD} and \cite{KLPDynamicsDomains}. Combining the results mentioned so far, we obtain: 

\begin{Thm}
  Let $\rho: F_g\rightarrow \PSL(2k,\bR)$ be a purely hyperbolic (Anosov) Schottky representation. Then, the properly discontinuous and cocompact action of $\rho$ on $D\subset \RP{2k-1}$ admits a fundamental domain bounded by finitely many polyhedral hypersurfaces.
\end{Thm}
Surprisingly, when $n=4k+3$ we can also build fundamental domains, but we have to pass to the double cover $S^{4k+2}$ of projective space in order to do so. The orbit $D\subset S^{4k+2}$  of such a fundamental domain coincides with the one predicted by the theory of domains of discontinuity in oriented flag manifolds recently developed in \cite{SteckerTreib}. We define halfspaces of the sphere $S^{4k+2}$ satisfying the same properties as the projective halfspaces above, and prove the following:
\begin{Thm}
  Let $\rho: F_g \rightarrow \PSL(4k+3,\bR)$ be a purely hyperbolic (Anosov) Schottky representation. Then, the properly discontinuous and cocompact action of $\rho$ on $D\subset S^{4k+2}$ admits a fundamental domain bounded by finitely many polyhedral hypersurfaces.
\end{Thm}
The main inspiration for defining halfspaces in spheres was the work of Choi and Goldman \cite{ChoiGoldmanTameness}. They use halfspaces to build fundamental domains in $S^2$ for Fuchsian representations in $\SO(2,1)$, and call the boundary of such a halfspace a \emph{crooked circle}. The quotients of $S^2$ obtained this way compactify quotients of $\bR^3$ by properly discontinuous affine actions of free groups. It might be possible to use cones over halfspaces in $S^{4k+2}$ in order to build fundamental domains for proper affine actions on $\bR^{4k+3}$, but this is outside the scope of this paper.

In order emphasize the similarities between positive and maximal representations, we have structured the paper in a similar way to the previous paper \cite{BTSchottky}.

We thank Daniele Alessandrini, Federica Fanoni, Misha Gekhtman, Fran\c cois Gu\'eritaud, Fanny Kassel, Giuseppe Martone, Beatrice Pozzetti, Anna-Sofie Schilling, Ilia Smilga, Florian Stecker, Anna Wienhard and Feng Zhu for insightful comments and helpful discussions.

\section{Preliminaries}
\subsection{Partially cyclically ordered spaces}    \label{sec:PCOSpaces}
In this section we recall definitions from \cite{BTSchottky} involving cyclic orders and the definition of a Schottky group in a cyclically ordered space.

\begin{Def}
A \emph{partial cyclic order} (PCO) on a set $C$ is a relation $\cycle{}{}{}$ on triples in $C$ satisfying, for any $a,b,c,d\in C$ :
\begin{itemize}
 \item if $\cycle{a}{b}{c}$, then $\cycle{b}{c}{a}$ (\emph{cyclicity});
 \item if $\cycle{a}{b}{c}$, then not $\cycle{c}{b}{a}$ (\emph{asymmetry});
 \item if $\cycle{a}{b}{c}$ and $\cycle{a}{c}{d}$, then $\cycle{a}{b}{d}$ (\emph{transitivity}).
\end{itemize}
If in addition the relation satisfies the following, then we call it a \emph{total cyclic order}:
\begin{itemize}
 \item If $a,b,c$ are distinct, then either $\cycle{a}{b}{c}$ or $\cycle{c}{b}{a}$ (\emph{totality}).
\end{itemize}
\end{Def}

\begin{Def}
A map $f : C \rightarrow D$ between partially cyclically ordered spaces $C,D$ is called \emph{increasing} if $\cycle{a}{b}{c}$ implies $\cycle{f(a)}{f(b)}{f(c)}$.
An automorphism of a partial cyclic order is an increasing map $f:C\rightarrow C$ with an increasing inverse. We will denote by $G$ the group of all automorphisms of $C$.
\end{Def} 

Partial cyclic orders give rise to a notion of intervals in $C$. Schottky groups in cyclically ordered spaces are modeled on the the case of Fuchsian Schottky groups acting on $\RP{1}$, and will be defined analogously using the following notion of interval.
\begin{Def}
Let $a,b\in C$. The \emph{interval} between $a$ and $b$ is the set 
\[ \ival{a}{b}:=\{x \in C ~|~ \cycle{a}{x}{b}\}. \]
The \emph{opposite} of an interval $I=\ival{a}{b}$ is the interval $\ival{b}{a}$, also denoted by $-I$.

The intervals in $C$ generate a natural topology under which order-preserving maps are continuous.
\end{Def}

\begin{Def}
We call a sequence $(a_n) \in C^\bN$ \emph{increasing} if $\cycle{a_i}{a_j}{a_k}$ whenever $i<j<k$.
\end{Def}

The cyclic order being only partial means that not every pair $a\neq b\in C$ is \emph{comparable}.

\begin{Def}\label{Def:comparableSet}
The \emph{comparable set} of a point $a\in C$ is
\[\comp(a) = \{x\in C ~|~ \ival{a}{x}\neq \emptyset \text{ or } \ival{x}{a}\neq \emptyset\}.\]
\end{Def}

Let $C$ be a partially cyclically ordered set. The following notion of completeness will ensure that Schottky groups defined using intervals have well defined limit sets.
\begin{Def}
$C$ is \emph{increasing-complete} if every increasing sequence converges to a unique limit in the interval topology.
\end{Def}

\begin{Def}
$C$ is \emph{proper} if for any increasing quadruple $(a,b,c,d)\in C^4$, we have $\overline{\ival{b}{c}} \subset \ival{a}{d}$. Here, ``bar'' denotes the closure in the interval topology.
\end{Def}


Let $\Gamma\subset \PSL(2,\bR)$ be a Schottky group acting on $\RP{1}$. That is, $\Gamma=\langle A_1,\dots,A_g\rangle$ where $A_j\in \PSL(2,\bR)$ and there exist $2g$ pairwise disjoint intervals $I_1^\pm,\dots,I_g^\pm\subset \RP{1}$ such that $A_j(-I_j^-)=I_j^+$. We will use $\Gamma$ as a combinatorial model for Schottky groups in general cyclically ordered spaces. Note that adjacent intervals are allowed to share an endpoint.

The \emph{limit set} of $\Gamma$ is the set of accumulation points of a $\Gamma$-orbit in $\RP{1}$. There are two possibilities for its topological type: it can be a Cantor set, or the whole projective line $\RP{1}$. 

In the latter case, we will say that $\Gamma$ is a \emph{finite area} model. This terminology comes from the fact that in this case the quotient of the hyperbolic plane $\Gamma\backslash\bH^2$ is a (noncompact) finite area hyperbolic surface. In a finite area model, each endpoint of a defining interval is necessarily shared with another interval. This setting will be used for the connection to positive representations (\autoref{sec:positive_reps}).

Another special case is when the defining intervals $I_j^\pm$ have disjoint closures. Any Schottky group $\Gamma$ admitting a presentation using intervals with disjoint closures is called \emph{purely hyperbolic}, because in this case every non-trivial element of the group is hyperbolic in $\PSL(2,\bR)$. Note that such a group always admits Schottky presentations where some (or all) endpoints are shared between two intervals as well. This setting will be used for the connection to Anosov representations and the construction of fundamental domains (\autoref{sec:Anosov_reps} and \autoref{sec:fundamentalD}).

There can be intermediate cases as well: If some ends of the surface $\Gamma\backslash\bH^2$ are cusps and some are funnels, the limit set is a Cantor set, but $\Gamma$ does not admit a presentation as a purely hyperbolic Schottky group. Apart from the general setup, we will not study such intermediate cases in this paper.

Let $G=\Aut(C)$ be the group of order-preserving bijections of $C$ with an order-preserving inverse.

\begin{Def} \label{Def:Gen_Schottky}
Let $\xi_0$ be an increasing map from the set of endpoints of the intervals $I_1^\pm,\dots,I_g^\pm$ into a partially cyclically ordered set $C$, where $I_i^\pm=\ival{a_i^\pm}{b_i^\pm}$. Define the corresponding image intervals in $C$ by $J_i^\pm=\ival{\xi_0(a_i^\pm)}{\xi_0(b_i^\pm)}$. Next, assume there exist $h_1,\dots,h_g\in G$ which pair the endpoints of $J_i^\pm$ in the same way that the $A_i$ pair the endpoints of $I^\pm_i$, so that $h_i(-J_i^-)=J_i^+$. We call the induced morphism $\rho:\Gamma \to G$ sending $A_i$ to $h_i$ a \emph{generalized Schottky representation}, its image in $G$ a \emph{generalized Schottky group} and the intervals $J_i^\pm$ used to define it a set of \emph{Schottky intervals} for this group.
\end{Def}

As a consequence of the \emph{Ping-Pong Lemma}, the generalized Schottky group $\rho(\Gamma)$ is freely generated by $h_1,\dots,h_g$. We can accurately describe the dynamics of the free group action using images of the defining intervals by group elements. We will define a bijection between length $k$ words in the group and certain intervals in the $\RP{1}$ model and in $C$.

For convenience, denote $A_{-j} = A_j^{-1}$ and $I_{-j}^+ = I_j^-$ for $j>0$. Let $\gamma=A_{j_1}\dots A_{j_k}$ be a reduced word of length $k$ in the generators $A_j$ and their inverses, where $j_i \in \{\pm 1,\ldots,\pm g\}$. Define $I_\gamma  = A_{j_1}\dots A_{j_{k-1}} I_{j_k}^+$.
We will call such an interval a $k$-th order interval. For example, first order intervals are associated to length $1$ words, which are just generators, and correspond to the defining Schottky intervals of the model.

Using the same construction with the intervals $J_j\subset C$ and generators $h_j=\rho(A_j)$, we define the $k$-th order interval $J_\gamma$ in $C$.

This bijection between words of length $k$ and $k$-th order intervals has the following property which will be useful when investigating infinite words:

\begin{Lem}\label{Lem:PCO_kthorder_nested}
If $\gamma\in \Gamma$ is a length $k$ reduced word and $\gamma=\gamma' A_i$ with $\gamma'$ a word of length $k-1$, then $J_\gamma \subset J_{\gamma'}$.
\begin{Prf}
Since the map $\xi_0$ is increasing, for any $k\neq -j$, we have $J_k^+ \subset -J_j^-$, thus $\rho(A_j) J_k^+ \subset J_j^+$. This means that if we denote the last letter of $\gamma'$ by $A_l$, 
\[J_\gamma = \rho(\gamma') J_i^+ = \rho(\gamma'' A_l) J_i^+ \subset \rho(\gamma'') J_l^+ = J_{\gamma'}. \qedhere \]
\end{Prf}
\end{Lem}

In \cite{BTSchottky}, the main theorem is the existence of an equivariant boundary map for a certain class of generalized Schottky groups :

\begin{Thm}\label{Thm:SchottkyLeftContinuousMap}
Let $\Gamma$ be a finite area model and $\rho: \Gamma \to G$ be a generalized Schottky representation. Assume that $C$ is first countable, increasing-complete and proper. Then there is a left-continuous, equivariant, increasing boundary map $\xi: \RP{1} \to C$.
\end{Thm}

This theorem provides a way to relate generalized Schottky groups to other interesting classes of representations which are defined by the existence of an equivariant boundary map.

In what follows, we introduce the space of oriented flags in $\bR^n$ and show that it admits a natural partial cyclic order invariant under $\PSL(n,\bR)$.
\subsection{Complete oriented flags}
We consider the vector space $\bR^n$, together with its standard basis and the induced orientation. Moreover, let $G=\PSL(n,\bR)$ and $B\subset G$ be the subgroup of upper triangular matrices.
\begin{Def}
A \emph{complete flag} $F$ in $\bR^n$ is a collection of nested subspaces
\[ \{0\} \subset F^{(1)} \subset F^{(2)} \subset \ldots \subset F^{(n-1)} \subset \bR^n, \]
where $\dim(F^{(i)}) = i$. For ease of notation, we sometimes include $F^{(0)} = \{0\}$ and $F^{(n)} = \bR^n$. We denote the space of complete flags by $\Flag{n}$.
\end{Def}
The group $G$ acts transitively on the space of complete flags. The stabilizer of the standard flag
\begin{equation}    \label{eq:std_flag}
    \langle e_1 \rangle \ \subset\  \langle e_1,e_2 \rangle \ \subset \ \ldots \ \subset \ \langle e_1,\ldots,e_{n-1}\rangle
\end{equation} 
is $B$, so the space of complete flags identifies with the homogeneous space $G/B$. We shall be interested in oriented flags. In terms of homogeneous spaces, this means that we consider the space $G/B_0$. Since $-1\in B_0$ if and only if $n$ is even, the space of complete oriented flags is a bit harder to describe in those dimensions. As an auxiliary object, we also consider the corresponding homogeneous space for the group $\SL(n,\bR)$.
\begin{Def}
  \begin{enumerate}
      \item A \emph{complete oriented flag for $\SL(n,\bR)$} is a complete flag in $\bR^n$ together with a choice of orientation on each of the subspaces $F^{(i)}, \ 1\leq i \leq n-1$. The space of complete oriented flags for $\SL(n,\bR)$ will be denoted $\loFlag{n}$.
      \item A \emph{complete oriented flag for $\PSL(n,\bR)$} is a complete flag in $\bR^n$ together with a choice of orientation on each of the subspaces $F^{(i)}, \ 1\leq i \leq n-1$, up to simultaneously reversing all the odd-dimensional orientations if $n$ is even. The space of complete oriented flags for $\PSL(n,\bR)$ will be denoted $\oFlag{n}$.
  \end{enumerate}
  The extremal dimension $F^{(n)}=\bR^n$ is always equipped with its standard orientation.
\end{Def}
\begin{Rem}
 Intuitively, the space $\loFlag{n}$ appears easier to describe and work with than $\oFlag{n}$. Moreover, in order to prove results about $\oFlag{n}$, we will frequently use lifts to $\loFlag{n}$. Our reason for using $\oFlag{n}$ is the better behavior of the partial cyclic order we are going to define. See \autoref{rem:reason_for_PSL} for more details on the problems that arise when using $\loFlag{n}$.
 
\end{Rem}
Again, $G$ acts transitively on $\oFlag{n}$. We can lift the standard flag \eqref{eq:std_flag} to $\oFlag{n}$ by equipping the $i$-dimensional component with the orientation determined by the ordered basis $(e_1,\ldots,e_i)$. Its stabilizer is $B_0$, yielding the identification
\[ \oFlag{n} = G/B_0. \]
The natural map $G \to G/B_0$ sends any element $g \in G$ to the image of the standard flag under $g$. In other words, $gB_0 \in G/B_0$ is the complete oriented flag $F_g$ such that the first $i$ columns of $g$ form an oriented basis for $F_g^{(i)}$ (up to simultaneously changing all odd-dimensional orientations if $n$ is even).

We will use matrices to denote elements of $\PSL(n,\bR)$ even though they are technically equivalence classes comprising two matrices if $n$ is even. Accordingly, statements such as ``all diagonal entries are positive'' should be interpreted as ``all diagonal entries are positive or all diagonal entries are negative''.
\subsection{Oriented transversality} \label{sec:or_transverse}
The first notion we require before we can define the partial cyclic order on oriented flags is an oriented version of transversality for flags. This notion appears under the name \emph{2-hyperconvexity} in \cite{guichard} and in the unoriented setting in \cite{labourie},\cite{GuichardHitchinHyperconvex}. We will need direct sums of oriented subspaces, so we first fix some notation.
\begin{Def}
  Let $V,W \subset \bR^n$ be oriented subspaces.
  \begin{itemize}
      \item If $V,W$ agree as oriented subspaces, we write $V \equalo W$;
      \item $-V$ denotes the same subspace with the opposite orientation;
      \item If $V$ and $W$ are transverse, we interpret $V \oplus W$ as an oriented subspace by equipping it with the orientation induced by the concatenation of a positive basis of $V$ and a positive basis of $W$, in that order.
  \end{itemize}
\end{Def}
Note that oriented direct sums depend on the ordering of the summands:
\[ V \oplus W \equalo (-1)^{\dim(V)\dim(W)} W\oplus V \]

\begin{Rem}
  We use negation to denote transformations of different spaces: On a fixed oriented Grassmannian, it denotes the involution inverting orientations. On the space $\loFlag{n}$ however, for even $n$, it denotes the induced action of $-1$ which inverts all odd-dimensional orientations.
\end{Rem}
\begin{Def} \label{Def:or_transverse}
Let $F_1,F_2\in\oFlag{n}$ be complete oriented flags. 
\begin{itemize}
    \item If $n$ is odd, the pair $(F_1,F_2)$ is called \emph{oriented transverse} if, for every $1\leq i\leq n-1$, we have
    \[ F_1^{(i)} \oplus F_2^{(n-i)} \equalo \bR^n; \]
    \item If $n$ is even, the pair $(F_1,F_2)$ is called \emph{oriented transverse} if there exist lifts $\widehat{F}_1,\widehat{F}_2 \in \loFlag{n}$ such that for every $1\leq i\leq n-1$, we have
    \[ \widehat{F}_1^{(i)} \oplus \widehat{F}_2^{(n-i)} \equalo \bR^n; \]
    We then call the pair $(\widehat{F}_1,\widehat{F}_2)$ a \emph{consistently oriented lift} of $(F_1,F_2)$.
    \item The set of flags that are oriented transverse to $F_1$ will be denoted by
    \[ \comp(F_1) = \{ F\in\oFlag{n} \mid (F_1,F) \ \text{is an oriented transverse pair} \} \]
    and, anticipative of the partial cyclic order, will be called the \emph{comparable set of $F_1$}.
\end{itemize}
\end{Def}
The left action of $\PSL(n,\bR)$ preserves oriented transversality. This is clear when $n$ is odd, and for even $n$ we just observe that $-1$ preserves the orientation of $\bR^n$. We also note that if $n$ is even and $F_1,F_2\in\oFlag{n}$ is an oriented transverse pair, there are exactly two consistently oriented lifts: If $(\widehat{F}_1,\widehat{F}_2)$ is one such pair, $(-\widehat{F}_1,-\widehat{F}_2)$ is the other.
\begin{Lem} \label{Lem:or_transv_symmetric}
Oriented transversality of pairs in $\oFlag{n}$ is symmetric.
\begin{Prf}
Let $(F_1,F_2)$ be an oriented transverse pair in $\oFlag{n}$, and let $(\widehat{F}_1,\widehat{F}_2)$ be a consistently oriented lift to $\loFlag{n}$ (if $n$ is odd, $\widehat{F}_i=F_i$). For each $i$, we have
\[ \widehat{F}_1^{(i)} \oplus \widehat{F}_2^{(n-i)} \equalo \bR^n. \]
It follows that
\[ \widehat{F}_2^{(n-i)} \oplus \widehat{F}_1^{(i)} \equalo (-1)^{i(n-i)} \bR^n. \]
If $n$ is odd, $(-1)^{i(n-i)}=1$, so $(F_2,F_1)$ is oriented transverse. If $n$ is even, $(-1)^{i(n-i)}=(-1)^i$. Then, the lifts $-\widehat{F}_2,\widehat{F}_1$ are consistently oriented and so $(F_2,F_1)$ is oriented transverse.
\end{Prf}
\end{Lem}
An example of an oriented transverse pair of flags, which we will call the standard pair, is given by $F_e = eB_0$, the identity coset, and
\[ F_{w_0} = w_0B_0 = \begin{pmatrix} & & & & \reflectbox{$\ddots$} \\ & & & -1 \\ & & 1 \\ & -1 \\ 1 \end{pmatrix}B_0.\]
To see that $(F_e,F_{w_0})$ is indeed an oriented transverse pair, observe that oriented transversality to $F_e$ is equivalent to all minors of $w_0$ obtained using the last $k$ rows and the first $k$ columns being positive. 

It will be very useful later on to choose special representatives for oriented flags. For brevity, we write
\[ \ltu = \left\{ g \in \PSL(n,\bR) \mid g \ \text{lower triangular and unipotent}  \right\} \]
and
\[ \utu = \left\{ g \in \PSL(n,\bR) \mid g \ \text{upper triangular and unipotent}  \right\}. \]
\begin{Lem} \label{Lem:representatives}
  Let $F\in\oFlag{n}$ such that $(F_{w_0},F)$ is an oriented transverse pair. Then $F$ admits a unique representative in $\ltu$. More precisely, the projection $G\to G/B_0 = \oFlag{n}$ restricts to a diffeomorphism 
  \[ \ltu \xrightarrow{\hspace{0.3cm}\cong\hspace{0.3cm}} \comp(F_{w_0}) \]
  with the comparable set of $F_{w_0}$.
  \begin{Prf}
  Let $F=gB_0$. Since the unoriented flag $gB \in \Flag{n}$ is transverse to $w_0B$, it is easy to see that $g$ can be chosen to be lower triangular. Oriented transversality ensures that all diagonal entries must be positive. If $gb$ is again lower triangular for some $b\in B_0$, $b$ is necessarily diagonal. Therefore, requiring $g$ to be unipotent fixes the representative uniquely, and the projection induces an injective map $f\colon\ltu\to G/B_0$ whose image is the (open) set of flags oriented transverse to $F_{w_0}$. The differential $\mathrm d_1f \colon\mathrm{Lie}(\ltu) \to \mathrm T_{[1]} G/B_0$ is an isomorphism, and by $G$-equivariance of the projection, $f$ is a diffeomorphism onto its image.
  \end{Prf}
\end{Lem}
\begin{Lem} \label{Lem:transitive_on_pairs}
  The left action of $\PSL(n,\bR)$ on oriented transverse pairs in $\oFlag{n}$ is transitive. The stabilizer of $(F_e,F_{w_0})$ is given by the subgroup $A$ of diagonal matrices with positive entries.
  \begin{Prf}
  Let $(F_1,F_2)$ be an oriented transverse pair in $\oFlag{n}$. Since the action of $\PSL(n,\bR)$ on $\oFlag{n}$ is transitive, we may assume that $F_2 = F_{w_0}$. Then, by \autoref{Lem:representatives}, $F_1 = gB_0$ for a unique representative $g\in\ltu$.
  The stabilizer of $F_{w_0}$ under the left action of $\PSL(n,\bR)$ is $B_0^t$ since $w_0B_0 = B_0^t w_0$.
  In particular, it contains the element $g^{-1}$ mapping $F_1$ to $F_e$. Since the stabilizer of $F_e$ under left multiplication is $B_0$, the stabilizer of the pair $(F_e,F_{w_0})$ is $A$, as claimed.
  \end{Prf}
\end{Lem}
\begin{Cor}\label{Cor:UniqueTransverseLift}
  Let $F_1,F_2\in\Flag{n}$ be a pair of transverse flags. Let $\hat{F}_1\in\oFlag{n}$ be a lift of $F_1$ to oriented flags. Then, there is a unique lift $\hat{F}_2\in\oFlag{n}$ of $F_2$ such that the pair $(F_1,F_2)$ is oriented transverse.
\end{Cor}
Since our description of oriented transversality in even dimension is based on choosing lifts to $\loFlag{n}$, it will be useful to describe some basic properties of oriented transversality in $\loFlag{n}$.

\begin{Def}
  Let $\widehat{F}_1,\widehat{F}_2 \in \loFlag{n}$ be complete oriented flags for $\SL(n,\bR)$. The pair $(\widehat{F}_1,\widehat{F}_2)$ is called \emph{oriented transverse} if we have
  \[ \widehat{F}_1^{(i)} \oplus \widehat{F}_2^{(n-i)} \equalo \bR^n \]
  for all $1\leq i\leq n-1$.
\end{Def}
Note that the identity matrix and $w_0$, considered as representatives of elements of $\loFlag{n}$, are oriented transverse. They will be our standard oriented transverse pair in $\loFlag{n}$. The following lemma shows how symmetry of oriented transversality fails in $\loFlag{n}$.
\begin{Lem} \label{Lem:or_transverse_SL_not_symmetric}
  Let $n$ be even. If $(\widehat{F}_1,\widehat{F}_2)$ is an oriented transverse pair in $\loFlag{n}$, then $(-\widehat{F}_2,\widehat{F}_1)$ is oriented transverse, and $(-\widehat{F}_1,\widehat{F}_2)$ is not.
  \begin{Prf}
  For each $i$, we have
  \[ \widehat{F}_2^{(n-i)} \oplus \widehat{F}_1^{(i)} \equalo (-1)^{i(n-i)} \widehat{F}_1^{(i)} \oplus \widehat{F}_2^{(n-i)} \equalo (-1)^{i(n-i)} \bR^n. \]
  The sign is negative if and only if $i$ is odd. Therefore, $(-\widehat{F}_2,\widehat{F}_1)$ is oriented transverse.\\
  To see that $(-\widehat{F}_1,\widehat{F}_2)$ is not oriented transverse, consider any splitting $\widehat{F}_1^{(i)} \oplus \widehat{F}_2^{(n-i)}$ where $i$ is odd.
  \end{Prf}
\end{Lem}
Let $\mathcal{E}=(e_1,\dots,e_n)$ be a positive ordered basis of $\bR^n$ (that is, it agrees with the standard orientation on $\bR^n$). We can associate a unique pair of oriented transverse flags $F_\mathcal{E}^+,F_\mathcal{E}^-$ to $\mathcal{E}$ in the following way:
\[F_\mathcal{E}^{+(k)} := \Span(e_1,\dots,e_k),\]
\[F_\mathcal{E}^{-(k)} := \Span(e_{n-k+1},\dots,e_n).\]
Here, each span is understood to be equipped with the orientation given by the ordering of basis vectors. Conversely, given a pair of oriented transverse flags $(F_1,F_2)$ we can find a positive ordered basis $\mathcal{E}$, unique up to multiplying each basis vector by a positive scalar, such that $F_1=F_\mathcal{E}^{+}$ and $F_2=F_\mathcal{E}^-$. We will say that such a basis is \emph{adapted to $F_1,F_2$}.\\
We can also associate the group $\ltu_\mathcal{E} = g\ltu g^{-1}$ to $\mathcal E$, where $g\in\PSL(n,\bR)$ is the change of basis from the standard basis to $\mathcal{E}$. Analogously to \autoref{Lem:representatives}, $\ltu_\mathcal{E}$ parametrizes the set of oriented flags oriented transverse to $F_2$.
%
\begin{Rem}
  Oriented transversality, as treated in this paper, is a special case of an \emph{oriented relative position}. These are all the possible combinatorial positions two oriented flags can be in -- more formally, an oriented relative position is a point in the quotient
  \[ \PSL(n,\bR) \backslash \left(\oFlag{n} \times \oFlag{n}\right) \]
  by the diagonal left-action of $\PSL(n,\bR)$. See \cite{SteckerTreib} for a more thorough treatment.
\end{Rem}

%
%
%

%
%
%
\subsection{Total positivity}
Total positivity of a matrix is a classical notion which has many applications (see e.g. \cite{Lusztig_survey}). In particular, it is used by Fock and Goncharov to define the higher Teichm\"uller spaces of positive representations. In this section, we recall some of the properties of total positivity.
\begin{Def}
Let $M\in\mathrm{Mat}(n,\bR)$ be a $(n\times n)$-matrix. Then $M$ is \emph{totally positive} if all minors of $M$ are positive.\\
If $M$ is either upper or lower triangular, we will call $M$ \emph{(triangular) totally positive} if all minors that do not vanish by triangularity are positive. Explicitly, if $M$ is upper (resp. lower) triangular, the minors to consider are determined by indices $i_1,\ldots i_k, \ j_1,\ldots j_k$ such that $i_l \leq j_l \ \forall l$ (resp. $i_l \geq j_l$).\\
An element of $\PSL(n,\bR)$ is called totally positive if it has a lift to $\SL(n,\bR)$ which is totally positive.\\
We will use the term \emph{totally nonnegative} in each of the cases above to denote the analogous situations where we only ask for minors to be nonnegative.
\end{Def}

We now introduce some notation for multiindices which will make the statement of many formulas involving minors simpler and more readable.
\begin{Def}
Let $k,n \in \bN$ be two integers. Then we write
\[ \mind{k}{n} := \{ (i_1,\ldots,i_k) \ | \ 1\leq i_1 < \ldots < i_k \leq n \} \]
for the set of multiindices with $k$ entries in increasing order from $\{1,\ldots,n\}$.
\end{Def}
Elements $\mathbf{i}\in\mind{k}{n}$ will be used to denote the rows or columns determining a minor: In combination with our previous notation, we can now write
\[ \submat{M}{\mathbf{i}}{\mathbf{j}} = \submat{M}{i_1\ldots i_k}{j_1\ldots j_k} \]
for the submatrix consisting of rows $i_1,\ldots,i_k$ and columns $j_1,\ldots,j_k$, and
\[ \minor{M}{\mathbf{i}}{\mathbf{j}} = \minor{M}{i_1\ldots i_k}{j_1\ldots j_k} = \det\left( \submat{M}{\mathbf{i}}{\mathbf{j}} \right) \]
for the corresponding minor.\\
The reason for using (ordered) multiindices instead of (unordered) k-subsets is that it makes them easier to compare.
\begin{Def}
Let $\mathbf{i} = (i_1,\ldots,i_k), \ \mathbf{j}=(j_1,\ldots,j_k) \in\mind{k}{n}$. We define a partial order $\leq$ on $\mind{k}{n}$ by
\[ \mathbf{i} \leq \mathbf{j} \Leftrightarrow i_l \leq j_l \ \forall l. \]
The \emph{absolute value} of a multiindex is the sum of its components,
\[ |\mathbf{i}| = \sum\limits_l i_l. \]
\end{Def}
The partial order on multiindices is particularly useful when working with triangular matrices. As mentioned earlier, if a matrix $M$ is upper (resp. lower) triangular, then all minors $\minor{M}{\mathbf{i}}{\mathbf{j}}$ with $\mathbf{i}> \mathbf{j}$ (resp. $\mathbf{i}< \mathbf{j}$) vanish automatically, and we call $M$ totally positive if $\minor{M}{\mathbf{i}}{\mathbf{j}} > 0 \ \forall \mathbf{i}\leq \mathbf{j}$ (resp. $\mathbf{i}\geq \mathbf{j}$).\\
As a first example of this notation in use, let us state the Cauchy-Binet formula. It describes how to calculate the determinant of a product of non-square matrices in terms of the minors of these matrices, and will play a central role later on. Using multiindices emphasizes the formal similarity to ordinary matrix multiplication (see for example \cite[(3.14)]{Tao} for a proof).
\begin{Lem}[Cauchy-Binet]
Let $M$ be a $(m\times r)$-matrix and $N$ a $(r\times m)$-matrix. Then, we have
\[ \det(MN) = \sum\limits_{\mathbf{k}\in\mind{m}{r}} \minor{M}{1\ldots m}{\mathbf{k}}\minor{N}{\mathbf{k}}{1\ldots m}. \]
\end{Lem}
Note that the formula includes the case $m>r$. Then $\det(MN)$ vanishes and, since $\mind{m}{r}$ is empty, the (empty) sum equals $0$ as well.\\
Our use of the formula lies in the calculation of minors of the product of two matrices: If $k\leq m$ and $\mathbf{i},\mathbf{j}\in\mind{k}{m}$ are multiindices, we obtain
\begin{equation}    \label{Eq:minor_product}
    \minor{(MN)}{\mathbf{i}}{\mathbf{j}} = \sum\limits_{\mathbf{k}\in\mind{k}{m}} \minor{M}{\mathbf{i}}{\mathbf{k}} \minor{N}{\mathbf{k}}{\mathbf{j}}.
\end{equation} 

As an immediate consequence of the Cauchy-Binet formula, one obtains the well-known fact that totally positive matrices form a semigroup.
\begin{Lem} \label{lem:totally_positive_semigroup}
Let $M,N\in\mathrm{Mat}(n,\bR)$ be totally positive. Then $MN$ is totally positive as well. If both $M$ and $N$ are upper (resp. lower) triangular and totally positive, then $MN$ is upper (resp. lower) triangular and totally positive.
\begin{Prf}
Let $M,N$ be totally positive. Then equation \eqref{Eq:minor_product} expresses any minor $\minor{(MN)}{\mathbf{i}}{\mathbf{j}}$ as a sum of positive summands.\\
If $M,N$ are both upper (resp. lower) triangular, then the same is true for $MN$, and for any two multiindices $\mathbf{i}\leq\mathbf{j} \in \mind{k}{n}$ (resp. $\mathbf{i}\geq\mathbf{j}$), we have
\[ \minor{(MN)}{\mathbf{i}}{\mathbf{j}} = \sum\limits_{\mathbf{k}\in\mind{k}{n}} \minor{M}{\mathbf{i}}{\mathbf{k}}\minor{N}{\mathbf{k}}{\mathbf{j}} = \sum\limits_{\substack{\mathbf{i}\leq\mathbf{k}\leq\mathbf{j} \\ \text{or}\ \mathbf{i}\geq\mathbf{k}\geq\mathbf{j}}} \minor{M}{\mathbf{i}}{\mathbf{k}}\minor{N}{\mathbf{k}}{\mathbf{j}}. \]
Since this sum is not empty, $MN$ is totally positive as well.
\end{Prf}
\end{Lem}

In \autoref{sec:fundamentalD} we will make use of the \emph{variation diminishing property} of totally positive matrices, introduced by Schoenberg \cite{schoenberg}.

\begin{Def} \label{def:variation}
  The \emph{upper (\emph{respectively} lower) variation} $\Var^+_\mathcal{E}(v)$ (resp. $\Var^-_\mathcal{E}(v)$) of a vector $x\in\bR^n$ with respect to an ordered basis $\mathcal{E}$ is the number of sign changes in the sequence of coordinates of $v$ in the basis $\mathcal{E}$, where $0$ coordinates are considered to have the sign which produces the largest value (respectively the lowest value).
  
  If we don't specify a basis $\mathcal{E}$, we mean the sign variation with respect to the canonical basis of $\bR^n$.
\end{Def}

\begin{Ex}
  $\Var^+(-1,2,0,3) = 3$ and $\Var^-(-1,2,0,3)=1$.
\end{Ex}


The variation diminishing property characterizes totally positive matrices in the following way.
\begin{Thm}[\cite{Pinkus}, Theorem 3.3]
\label{Thm:VariationDiminishing}
Let $A$ be an $n\times n$ totally positive matrix. Then, for any nonzero vector $v\in \bR^n$, we have
  \begin{itemize}
      \item $\Var^+(Av)\leq\Var^-(v)$;
      \item If $\Var^+(Av)=\Var^-(v)$, the sign of the last nonzero component of $v$ is the same as that of the last component of $Av$. If the last component of $Av$ is zero, the sign used when determining $\Var^+(Av)$ is used instead (see \autoref{rem:lastsignwelldefined}).
  \end{itemize}
    Conversely, any matrix $A$ with these properties is totally positive.
\end{Thm}
\begin{Rem}\label{rem:lastsignwelldefined}
  Let $v$ be a nonzero vector such that its $k$-th component $v_k$ is nonzero and $v_i$ vanishes for $i>k$. Then, there is a unique way of assigning signs to the last $n-k$ zeroes such that the variation is maximized. In particular, this gives the last component of $Av$ a well-defined sign.
\end{Rem}
The decomposition theorem for lower triangular totally positives matrices, due to A. Whitney \cite{WhitneyDecomp} and generalized by Lusztig \cite{lusztigtotpos}, will also be useful :
\begin{Thm}
  \label{Thm:TotPosDecomposition}
  Let $L$ be a unipotent, lower triangular $n\times n$ matrix. Then, $L$ is totally positive if and only if it can factored as
  \[C_1 \dots C_{n-2}C_{n-1}\]
  where
  \[C_r(\alpha_r^1,\dots,\alpha_r^{r-1}) = x_{r}(\alpha_r^{r})x_{r-1}(\alpha_r^{r-1})\dots x_1(\alpha_r^1).\]
  for some $\alpha_k^i >0$ and $x_k(\alpha)$ is an $n\times n$ matrix with $1$s on the diagonal, $\alpha$ in the entry $(k+1,k)$ and zero elsewhere.
\end{Thm}
\begin{Rem}
  In fact, in the previous theorem, the order in which the $x_k(\alpha)$ are multiplied can be chosen in many ways. More precisely, denote by $s_i$ the permutation $(i, i+1)$. Then, $s_1,\dots,s_{n-1}$ are generators for the symmetric group on $n$ letters. For any minimal expression $w=s_{i_1}\dots s_{i_L}$ of the longest word $w\in S_n$ in these generators, the corresponding product
  \[L = x_{i_1}(\alpha_{i_1})x_{i_2}(\alpha_{i_2})\dots x_{i_L}(\alpha_{i_L})\]
  is lower triangular totally positive whenever $\alpha_{i_k}>0$. In the theorem above we chose the expression
  \[w = s_1 (s_2 s_1) \dots (s_{n-1}\dots s_2 s_1).\]
\end{Rem}
%
%
%
%
%
%
%
%

\section{A partial cyclic order on oriented flags}   \label{Sec:PCO}
\subsection{Oriented 3-hyperconvexity}
The following property of triples of flags is the core of the partial cyclic order we are going to define. This is an oriented version of Fock-Goncharov's triple positivity \cite{FockGoncharov} and Labourie's $3$-hyperconvexity \cite{labourie}.
\begin{Def} \label{Def:hyperconvex}
Let $(F_1,F_2,F_3)$ be a triple in $\oFlag{n}$.
\begin{itemize}
    \item If $n$ is odd, $(F_1,F_2,F_3)$ is called \emph{oriented $3$-hyperconvex} if, for every triple of integers $0\leq i_1,i_2,i_3 \leq n-1$ satisfying $i_1+i_2+i_3 = n$,
    \[ F_1^{(i_1)} \oplus F_2^{(i_2)} \oplus F_3^{(i_3)} \equalo \bR^n; \]
    \item If $n$ is even, $(F_1,F_2,F_3)$ is called \emph{oriented $3$-hyperconvex} if there exist lifts $\widehat{F}_1,\widehat{F}_2,\widehat{F}_3 \in\loFlag{n}$ such that, for every triple of integers $0\leq i_1,i_2,i_3 \leq n-1$ satisfying $i_1+i_2+i_3 = n$;
    \[ \widehat{F}_1^{(i_1)} \oplus \widehat{F}_2^{(i_2)} \oplus \widehat{F}_3^{(i_3)} \equalo \bR^n. \]
    We call the triple $(\widehat{F}_1,\widehat{F}_2,\widehat{F}_3)$ a \emph{consistently oriented lift} of $(F_1,F_2,F_3)$.
\end{itemize}
\end{Def}
Since oriented 3-hyperconvexity for triples of oriented flags is the only notion of hyperconvexity appearing in this paper, we will simply call such triples \emph{hyperconvex}. Note that allowing one of the $i_j$ to vanish automatically includes oriented transversality of $(F_1,F_2)$, $(F_1,F_3)$ and $(F_2,F_3)$ in the definition. Like oriented transversality, hyperconvexity is invariant under the action of $\PSL(n,\bR)$ on $\oFlag{n}$. Moreover, in even dimension, if $(F_1,F_2,F_3)$ is a hyperconvex triple, it has exactly two consistently oriented lifts which are related by applying $-1$ to all its elements.

\autoref{Lem:representatives} showed that oriented flags in $\comp(F_{w_0})$ admit unique representatives in $\ltu$. We now examine the additional properties this representative satisfies if $(F_e,F,F_{w_0})$ is a hyperconvex triple.
%
%
We will use the notation
\[ \ltup = \{ g\in\ltu \mid g \ \text{totally positive} \}, \qquad \ltun = \{ g\in\ltu \mid g \ \text{totally nonnegative} \} \]
and
\[ \utup = \{ g\in\utu \mid g \ \text{totally positive} \}, \qquad \utun = \{ g\in\utu \mid g \ \text{totally nonnegative} \}. \]
%
%
%
%

\begin{Lem} \label{lem:std_triple}
  The projection $G\to G/B_0 = \oFlag{n}$ restricts to a diffeomorphism
  \[ \ltup \xrightarrow{\hspace{0.3cm}\cong\hspace{0.3cm}} \{ F\in\oFlag{n} \mid (F_e,F,F_{w_0}) \ \text{hyperconvex} \}. \]
  
  \begin{Prf}
  Using \autoref{Lem:representatives} and the fact that $\ltup$ is open in $\ltu$, the only thing left to show is that the preimage of the right hand side under the diffeomorphism $\ltu \xrightarrow{\cong} \comp(F_{w_0})$ is $\ltup$.
  
  Assume that $(F_e,F,F_{w_0})$ is hyperconvex. If $n$ is even, let $\lift{F}_e = eB_0^{\SL(n,\bR)}\in \loFlag{n}$, $\lift{F}_{w_0} = w_0B_0^{\SL(n,\bR)} $ and $\lift{F}\in\loFlag{n}$ the lift such that $(\lift{F}_e,\lift{F},\lift{F}_{w_0})$ is a consistently oriented lift (if $n$ is odd, $\lift{F}=F$). Let $M\in\SL(n,\bR)$ be any matrix representative for $\lift{F}$. Then the conditions on $M$ are as follows: Let $i_1,i_2,i_3$ be a triple of nonnegative integers satisfying $i_1+i_2+i_3 = n$. The oriented direct sum condition of \autoref{Def:hyperconvex} means that the matrix composed of the first $i_1$ columns of the identity, the first $i_2$ columns of $M$ and the first $i_3$ columns of $w_0$ (in that order) has positive determinant. We write $I_j$ for the $j\times j$ identity matrix and
\[ J_j = \begin{pmatrix} & & & \reflectbox{$\ddots$} \\ & & 1 \\ & -1 \\ 1 \end{pmatrix} \]
for the $j\times j$ antidiagonal matrix with alternating entries $\pm 1$, starting with $+1$ in the lower left corner. The matrix in question has the form
\[ \begin{pmatrix} I_{i_1} & \ast &  \\ & \submat{M}{i_1+1 \ldots i_1+i_2}{1\ldots i_2} &  \\  & \ast & J_{i_3} \end{pmatrix}, \]
where the stars are irrelevant for calculating the determinant. Since $J_j$ has determinant $1$, hyperconvexity of the triple is equivalent to 
\begin{equation}    \label{Eq:std_triple}
    \minor{M}{i_1+1 \ldots i_1+i_2}{1\ldots i_2} > 0 \qquad \forall i_1\geq 0, \ i_2 \geq 1, \ i_1+i_2 \leq n.
\end{equation} 
Now assume that $M$ is the unique representative in $\ltu$. By \cite[Theorem 2.8]{Pinkus}, positivity of all ``left--bound'' connected minors, that is, all the minors appearing in \eqref{Eq:std_triple}, is already sufficient to conclude that this representative is (triangular) totally positive.
  \end{Prf}
\end{Lem}

\begin{Cor}
The stabilizer in $\PSL(n,\bR)$ of a hyperconvex triple is trivial.
\begin{Prf}
By transitivity of the action on oriented transverse pairs (\autoref{Lem:transitive_on_pairs}), we can assume that the triple is of the form $(F_e,F,F_{w_0})$. Moreover, the stabilizer of $(F_e,F_{w_0})$ is $A$, the identity component of diagonal matrices. By \autoref{lem:std_triple}, $F=gB_0$ for a unique representative $g\in\ltup$. For any $a\in A$, the representative in $\ltup$ of the image $agB_0$ is given by $aga^{-1}$. Since conjugation by $a$ has a fixed point in $\ltup$ if and only if $a=1$, the claim follows.
\end{Prf}
\end{Cor}
A second corollary of \autoref{lem:std_triple} is the very simple relation between \emph{positivity} of triples of unoriented flags (from \cite{FockGoncharov}) and hyperconvexity of triples of oriented flags.
\begin{Def}\label{def:positivity}
A triple of unoriented flags, considered as cosets in $G/B$, is \emph{positive} if it is in the $G$-orbit of a triple of the form $B,gB,w_0B$ for some $g\in\ltup$.

More generally, an $N$-tuple of flags $F_1,\dots,F_N$ is positive if it is in the $G$-orbit of a tuple of the form
\[B,g_1B,g_1g_2B,\dots,g_1g_2\dots g_{N-2}B,w_0B,\]
with $g_i\in \ltup$.
\end{Def}
\begin{Cor}\label{cor:PositiveLiftHyperconvex}
A triple of unoriented flags if positive if and only if it admits a hyperconvex lift to oriented flags.
\end{Cor}
To close out this subsection, we compare the set of oriented flags $F\in\oFlag{n}$ such that $(F_e,F,F_{w_0})$ is hyperconvex with the intersection $\comp(F_e)\cap\comp(F_{w_0})$.
\begin{Prop}\label{prop:IntervalIsConnectedComponent}
The set $\{ F\in\oFlag{n} \mid (F_e,F,F_{w_0}) \ \text{hyperconvex} \}$ is a connected component of $\comp(F_e)\cap\comp(F_{w_0})$.
\begin{Prf}
Since hyperconvexity is an open condition, it suffices to prove that $\{ F\in\oFlag{n} \mid (F_e,F,F_{w_0}) \ \text{hyperconvex} \}$ is closed in $\comp(F_e)\cap\comp(F_{w_0})$. Let $F_t = g_tB_0, t \in [0,1],$ be a path in $\comp(F_e)\cap\comp(F_{w_0})$, with $g_t \in \ltu$. Assume that $(F_e,F_t,F_{w_0})$ is a hyperconvex triple for $t<1$, so $g_t \in \ltup$ for $t<1$. Then $g_1 \in \ltun$, and since $F_1 \in \comp(F_e)$, all ``bottom-left'' minors of $g_1$ are positive,
\[ \minor{(g_1)}{n-k+1\ldots n}{1\ldots k} > 0, \qquad 1\leq k \leq n \]
(as usual, for $n$ even, this holds for one of the two lifts to $\SL(n,\bR)$). Positivity of these minors and total nonnegativity of $g_1$ already implies $g_1 \in \ltup$ by \cite[Proposition 2.9]{Pinkus}.
\end{Prf}
\end{Prop}

\subsection{The partial cyclic order}

We now show that hyperconvexity satisfies the axioms of a partial cyclic order.
\begin{Prop}    \label{Prop:PCO_odd}
The relation $\mathcal{R} \subset (\oFlag{n})^3$ defined by
\[ (F_1,F_2,F_3) \in \mathcal{R} \quad \Leftrightarrow \quad (F_1,F_2,F_3) \ \text{is hyperconvex} \]
is a partial cyclic order. We will use the arrow notation from \autoref{sec:PCOSpaces} to denote it.
\begin{Prf}
Assume that $(F_1,F_2,F_3) \in \mathcal{R}$. If $n$ is even, let $(\widehat{F}_1,\widehat{F}_2,\widehat{F}_3)$ denote a consistently oriented lift.\\ 
We first check that the relation is asymmetric. Let us start with the odd-dimensional case, since it does not involve choices of lifts. For any $i_1,i_2,i_3$, we have
\[ F_1^{(i_1)} \oplus F_2^{(i_2)} \oplus F_3^{(i_3)} \equalo \bR^n \]
and therefore
\[ F_3^{(i_3)} \oplus F_2^{(i_2)} \oplus  F_1^{(i_1)} \equalo (-1)^{i_1(n-i_1) + i_2i_3} \bR^n. \]
Whenever $i_2$ and $i_3$ are both odd, we get the negative sign, showing that $(F_3,F_2,F_1) \not\in \mathcal{R}$.\\
If $n$ is even, we want to show that the triple $(F_3,F_2,F_1)$ does not admit a consistently oriented lift. Assume for the sake of contradiction that such a lift exists. Without loss of generality, it contains $\widehat{F}_1$. Then by \autoref{Lem:or_transverse_SL_not_symmetric}, the pairs $(-\widehat{F}_2,\widehat{F}_1)$ and $(-\widehat{F}_3,\widehat{F}_1)$ are oriented transverse, so the lift of the triple must contain $-\widehat{F}_2$ and $-\widehat{F}_3$. But another application of \autoref{Lem:or_transverse_SL_not_symmetric} shows that $(-\widehat{F}_3,-\widehat{F}_2)$ is not oriented transverse, a contradiction.

Now we turn to cyclicity. Again, we first treat the case of odd dimension. Let $i_1,i_2,i_3$ be integers such that $i_1+i_2+i_3 = n$. Then we have
\[ F_2^{(i_2)} \oplus F_3^{(i_3)} \oplus F_1^{(i_1)} \equalo (-1)^{i_1(i_2+i_3)} F_1^{(i_1)} \oplus F_2^{(i_2)} \oplus F_3^{(i_3)} \equalo (-1)^{i_1(n-i_1)} \bR^n. \]
As $n$ is odd, $i_1(n-i_1)$ is always even, and we conclude that $(F_2,F_3,F_1)\in\mathcal{R}$.\\
If $n$ is even, the same calculation with $\widehat{F}_i$ instead of $F_i$ yields a negative sign whenever $i_1$ is odd. Therefore, $(\widehat{F}_2,\widehat{F}_3,-\widehat{F}_1)$ is a consistently oriented lift of $(F_2,F_3,F_1)$.

Finally, we prove transitivity using the semigroup structure of totally positive matrices. Assume that we have a fourth flag $F_4\in\oFlag{n}$ such that $(F_1,F_3,F_4)$ is a hyperconvex triple. By \autoref{Lem:or_transv_symmetric} and \autoref{Lem:transitive_on_pairs}, we can use the $\PSL(n,\bR)$ action to normalize $F_1 = F_{w_0}$ and $F_2 = F_e$. Then, by oriented transversality with $F_{w_0}$, we have $F_3=g_3B_0, F_4=g_4B_0$ with representatives $g_3,g_4 \in \ltu$. Cyclicity implies that the triple $(F_e,F_3,F_{w_0})$ is hyperconvex, so \autoref{lem:std_triple} shows that $g_3$ is totally positive. Now consider the left-action of $g_3^{-1}$ on $\oFlag{n}$. It maps the triple $(F_{w_0},F_3,F_4)$ to $(F_{w_0},F_e,g_3^{-1}F_4)$. $g_3^{-1}F_4$ is represented by $g_3^{-1}g_4$, thus another application of cyclicity and \autoref{lem:std_triple} shows that $g_3^{-1}g_4$ is totally positive. Since totally positive matrices form a semigroup, we conclude that $g_4 = g_3(g_3^{-1}g_4)$ is totally positive as well. The triple $(F_{w_0},F_e,F_4) = (F_1,F_2,F_4)$ is therefore hyperconvex and the proof is complete.
\end{Prf}
\end{Prop}
%
\begin{Rem} \label{rem:reason_for_PSL}
  Let us compare how a similar construction for the space $\loFlag{n}$ behaves if $n$ is even. In \autoref{Def:hyperconvex}, a notion of hyperconvexity for triples in $\loFlag{n}$ appeared (though not explicitly named): Call a triple $\widehat F_1, \widehat F_2, \widehat F_3 \in \loFlag{n}$ hyperconvex if, for all $0\leq i_1,i_2,i_3 \leq n-1$ satisfying $i_1+i_2+i_3 = n$,
  \[ \widehat{F}_1^{(i_1)} \oplus \widehat{F}_2^{(i_2)} \oplus \widehat{F}_3^{(i_3)} \equalo \bR^n. \]
  On the space $\loFlag{n}$, this relation does not satisfy the cyclicity axiom of a partial cyclic order (see \autoref{Lem:or_transverse_SL_not_symmetric}).\\  
  One can circumvent this issue by defining the triple to be increasing if either of the triples $(\widehat F_1, \widehat F_2, \widehat F_3), (-\widehat F_1, \widehat F_2, \widehat F_3), (\widehat F_1, \widehat F_2, -\widehat F_3)$ is hyperconvex. This does indeed define a partial cyclic order on $\loFlag{n}$, but the somewhat artificial construction has an unwanted side effect: There are two distinct types of intervals. If the endpoints $(F_1,F_3)$ are oriented transverse, the interval $\ival{F_1}{F_3}$ is homeomorphic to its projection to $\oFlag{n}$. If $(F_1,-F_3)$ are oriented transverse, however, the interval decomposes into the disjoint union 
  \[ \ival{F_1}{F_3} = \ival{F_1}{-F_3} \sqcup \ival{-F_1}{F_3}. \]
  Note that the induced action of $-1\colon\bR^n\to\bR^n$ is a homeomorphism between the two intervals on the right hand side. An interval and its opposite are always of the two different types, and an interval of the connected type never contains a pair $(\widehat F,-\widehat F)$. This makes the construction of Schottky groups impossible for this partial cyclic order: The opposite of any Schottky interval needs to contain all the other Schottky intervals, but an interval of the connected type cannot contain an interval of the disconnected type.
\end{Rem}
%
%
In the proof of \autoref{Prop:PCO_odd}, we obtained the following useful characterization of cycles.
\begin{Lem} \label{Lem:transitivity}
Let $(F_e,F_1,F_2,F_{w_0})$ be a cycle, and let $F_1=g_1B_0,F_2=g_2B_0$, where $g_1,g_2 \in \ltup$. Then we have $g_2=g_1 h$ with $h \in \ltup$.
\end{Lem}
If $(F_1,F_3)$ is an oriented transverse pair, the interval between them is
\[ \ival{F_1}{F_3} = \{ F\in\oFlag{n} \ | \ \cycle{F_1}{F}{F_3} \}. \]
In \autoref{lem:std_triple}, we saw that the interval $\ival{F_e}{F_{w_0}}$ is given by the set $\ltup$ of all lower triangular, unipotent, totally positive matrices. It will be useful later on to have a similar description of the opposite interval $\ival{F_{w_0}}{F_e}$. In order to obtain this description, we first associate an involution $\chinv{F_1}{F_2}$ to any oriented transverse pair $(F_1,F_2)$. This involution will fix $F_1$ and $F_2$ and reverse the PCO, thereby providing a kind of symmetry for increasing and decreasing sequences.
%
\begin{Def}
Let $t\in\PSL(n,\bR)$ be the diagonal matrix with alternating $\pm 1$ entries,
\[ t = \begin{pmatrix} 1 \\ & -1 \\ & & 1 \\ & & & \ddots \end{pmatrix}. \]
The involution $\chinv{F_e}{F_{w_0}}$ is defined by
\begin{align*}
    \chinv{F_e}{F_{w_0}}\colon \oFlag{n} & \to \oFlag{n} \\
    gB_0 & \mapsto tgtB_0.
\end{align*}
If $(F_1,F_2)$ is an oriented transverse pair, pick an element $h\in\PSL(n,\bR)$ mapping $(F_e,F_{w_0})$ to $(F_1,F_2)$ and define
\[ \chinv{F_1}{F_2} = h\left(\chinv{F_e}{F_{w_0}}\right)h^{-1}. \] 
\end{Def}
First of all, we observe that $tB_0t = B_0$, so $\chinv{F_e}{F_{w_0}}$ is well-defined. Furthermore, if $d\in A = \mathrm{Stab}(F_e,F_{w_0})$,
\[ d\chinv{F_e}{F_{w_0}}d^{-1}(gB_0) = d t d^{-1} gtB_0 = tgtB_0. \]
Consequently, the involution $\chinv{F_1}{F_2}$ does not depend on the choice of $h$ in the definition above, but only on the pair $(F_1,F_2)$. Different involutions are related by 
\[ k\chinv{F_1}{F_2}k^{-1} = \chinv{kF_1}{kF_2} \] 
for any element $k\in\PSL(n,\bR)$.
\begin{Lem} \label{Lem:involution}
Let $(F_1,F_2)$ be an oriented transverse pair in $\oFlag{n}$ and $F\in\ival{F_1}{F_2}$. Then 
$\chinv{F_1}{F_2}(F)\in\ival{F_2}{F_1}$.
\begin{Prf}
We first show that it is enough to consider $(F_1,F_2)=(F_e,F_{w_0})$. Let $h\in\PSL(n,\bR)$ be such that $(hF_e,hF_{w_0})=(F_1,F_2)$. Then $h^{-1} \ival{F_1}{F_2} = \ival{F_e}{F_{w_0}}$, $\chinv{F_1}{F_2}=h\chinv{F_e}{F_{w_0}}h^{-1}$, and 
\[ \chinv{F_1}{F_2}(F)\in\ival{F_2}{F_1} \quad \Leftrightarrow \quad \chinv{F_e}{F_{w_0}}(h^{-1}F) \in \ival{F_{w_0}}{F_e}. \]
Next, let $F=gB_0$, where $g\in \ltup$, so that $\chinv{F_e}{F_{w_0}}F = tgtB_0$. By equivariance of the partial cyclic order, $tgtB_0\in \ival{F_{w_0}}{F_e}$ if and only if $\cycle{((tgt)^{-1}B_0)}{F_e}{F_{w_0}}$, which by cyclicity is equivalent to $tg^{-1}t\in \ival{F_{w_0}}{F_e}$. From \autoref{lem:std_triple}, this is equivalent to $t g^{-1} t\in \ltup$. The minors of the inverse of a matrix satisfy
\[(g^{-1})^\mathbf{i}_\mathbf{j} = (-1)^{|\mathbf{i}|+|\mathbf{j}|}\frac{g^{\mathbf{i}'}_{\mathbf{j}'}}{\det(g)},\]
where $\mathbf{i},\mathbf{i}'$ and $\mathbf{j},\mathbf{j}'$ are complementary multiindices, in the sense that $\mathbf{i}\sqcup\mathbf{i}'=\{1,\dots,n\}$ \cite[Section 1.1]{Pinkus}. Since conjugation by $t$ precisely multiplies a minor $g^\mathbf{i}_\mathbf{j}$ by $(-1)^{|\mathbf{i}|+|\mathbf{j}|}$, a matrix $g$ is in $\ltup$ if and only if $tg^{-1}t$ is in $\ltup$, proving the lemma.

\end{Prf}
\end{Lem}

\begin{Cor}
\label{Cor:standard_opposite_triple}
Let $F \in\oFlag{n}$ be a complete oriented flag such that $\left( F_{w_0} , F , F_e \right) $ is a hyperconvex triple. Then $F = gB_0$ with the unique representative $g \in \ltu$ satisfying
\[ (-1)^{|\mathbf{i}| + |\mathbf{j}|} \minor{g}{\mathbf{i}}{\mathbf{j}} > 0 \qquad \forall \mathbf{i}\geq\mathbf{j} \in\mind{k}{n} \ \forall k\leq n \]
(if $n$ is even, these equations are understood to hold for the unipotent lift to $\SL(n,\bR)$). Conversely, if $F$ has such a representative, the triple is hyperconvex.
\end{Cor}

\begin{Cor}
Let $(F_1,F_2)$ be an oriented transverse pair in $\oFlag{n}$. Then $\chinv{F_1}{F_2}$ reverses the partial cyclic order.
\begin{Prf}
Let $(F_3,F_4)$ be another oriented transverse pair. To improve readability, for the proof of this corollary, we set $\tau_{1,2} := \chinv{F_1}{F_2}$ and $\tau_{3,4} := \chinv{F_3}{F_4}$. We want to show that $\tau_{1,2} \ival{F_3}{F_4} = \ival{\tau_{1,2}F_4}{\tau_{1,2}F_3}$. By the previous lemma, we have
\[ \tau_{1,2} \ival{F_3}{F_4} = \tau_{1,2}\tau_{3,4} \ival{F_4}{F_3}. \]
The composition $\tau_{1,2}\tau_{3,4}$ is realized by the action of an element of $\PSL(n,\bR)$: If $h,k\in\PSL(n,\bR)$ are chosen such that $(hF_e,hF_{w_0})=(F_1,F_2)$ and $(kF_e,kF_{w_0})=(F_3,F_4)$, we know that
\[ \tau_{1,2}\tau_{3,4}(gB_0) = ktk^{-1}hth^{-1}gttB_0 = (ktk^{-1}hth^{-1})gB_0. \]
Therefore,
\[ \tau_{1,2}\tau_{3,4} \ival{F_4}{F_3} = \ival{\tau_{1,2}\tau_{3,4}F_4}{\tau_{1,2}\tau_{3,4}F_3} = \ival{\tau_{1,2}F_4}{\tau_{1,2}F_3}. \]
\end{Prf}
\end{Cor}

We now turn to proving a list of useful properties of the partial cyclic order on $\oFlag{n}$ given by hyperconvexity. These properties will allow us to construct Schottky representations and use the corresponding results from \cite{BTSchottky}.

\begin{Prop}
The partial cyclic order on $\oFlag{n}$ determined by hyperconvexity is increasing-complete and proper.
\begin{Prf}
We first prove properness. Using the action of $\PSL(n,\bR)$, we can bring an arbitrary $4$-cycle into the form $(F_e,F,F',F_{w_0})$. We want to show that $\overline{\ival{F}{F'}} \subset \ival{F_e}{F_{w_0}}$. Let $F=gB_0$ and $F'=g'B_0$ with $g,g' \in \ltup$ (see \autoref{lem:std_triple}). Let $F_m = g_mB_0 \in\ival{F}{F'}$ be a sequence converging to a flag $F_\infty$, where $g_m \in \ltup$. By transitivity of the PCO, we know that $(F_e,F,F_m,F_{w_0})$ and $(F_e,F_m,F',F_{w_0})$ are cycles. Therefore, \autoref{Lem:transitivity} shows that $g_m=gh_m$ and $g'=g_mh_m'$ for some $h_m,h_m' \in \ltup$. We denote the totally positive lifts to $\SL(n,\bR)$ by decorating with a hat. Then it is not hard to see that minors of $\widehat{g},\widehat{g}_n,\widehat{g}'$ are ordered the same way as the flags: For multiindices $\mathbf{i}\geq\mathbf{j}\in\mind{k}{n}$, we have
\[ \minor{\widehat{g}}{\mathbf{i}}{\mathbf{j}} \leq \minor{(\widehat{g}_n)}{\mathbf{i}}{\mathbf{j}} \leq \minor{(\widehat{g}')}{\mathbf{i}}{\mathbf{j}} \]
(see \autoref{Lem:minors_increasing} for a proof). In fact, strict inequality holds unless $\mathbf{i}=\mathbf{j}$, in which case all three minors are equal to $1$. Choosing singletons in the inequalities above bounds the matrix entries of $\widehat{g}_n$ between those of $\widehat{g}$ and $\widehat{g}'$. Since the sequence $F_n$ was assumed to converge to $F_\infty$, we conclude that the matrix entries of $\widehat{g}_n$ converge and $F_\infty$ is represented by the limit $\widehat{g}_\infty \in \ltu$. Furthermore, all minors $\minor{(\widehat{g}_\infty)}{\mathbf{i}}{\mathbf{j}}$ lie between the corresponding minors of $\widehat{g}$ and $\widehat{g}'$, so $\widehat{g}_\infty$ is totally positive. This shows that $F_\infty \in \ival{F_e}{F_{w_0}}$ and completes the proof of properness.\\
For increasing-completeness, assume that $F_1,F_2,\ldots$ is an increasing sequence. Let us normalize so that $F_2=F_{w_0}$ and $F_3=F_e$. Since we have $\cycle{F_{w_0}}{F_e}{F_m}$ for every $m\neq2,3$, we have $F_m = g_mB_0$ with $g_m \in \ltup$. Moreover, by \autoref{Lem:transitivity}, we have $g_{m+1}=g_mh_m$ for $m\geq 3$ and $g_1 = g_mh_{m,1}$, where $h_m, h_{m,1} \in \ltup$. Lifting to $\SL(n,\bR)$ as before, this implies that the entries below the diagonal of $\widehat{g}_m, \ m>4$ are bounded between the corresponding entries of $\widehat{g}_4$ and $\widehat{g}_1$, and they are strictly increasing in $m$. Therefore, there exists a unique limit $F_\infty = \lim_{m\to\infty}F_m$.
\end{Prf}
\end{Prop}

\begin{Lem} \label{Lem:minors_increasing}
Let $A,B \in \SL(n,\bR)$ be lower triangular, unipotent, totally positive. Let $\mathbf{i}\geq\mathbf{j}\in\mind{k}{n}$ be two multiindices. Then we have $\minor{(AB)}{\mathbf{i}}{\mathbf{j}}\geq \max(\minor{A}{\mathbf{i}}{\mathbf{j}},\minor{B}{\mathbf{i}}{\mathbf{j}})$. The inequality is strict unless $\mathbf{i}=\mathbf{j}$.
\begin{Prf}
The Cauchy-Binet formula yields
\[ \minor{(AB)}{\mathbf{i}}{\mathbf{j}} = \sum\limits_{\mathbf{k}\in\mind{k}{n}} \minor{A}{\mathbf{i}}{\mathbf{k}} \minor{B}{\mathbf{k}}{\mathbf{j}} = \sum\limits_{\mathbf{i}\geq\mathbf{k}\geq\mathbf{j}} \minor{A}{\mathbf{i}}{\mathbf{k}} \minor{B}{\mathbf{k}}{\mathbf{j}}. \]
All summands are positive, hence we obtain the following lower bound by only considering the two summands where $\mathbf{k}=\mathbf{i}$ or $\mathbf{k}=\mathbf{j}$:
\begin{align*}
    \text{If} \ \mathbf{i}=\mathbf{j}: & \quad \minor{(AB)}{\mathbf{i}}{\mathbf{i}} = \minor{A}{\mathbf{i}}{\mathbf{i}} = \minor{B}{\mathbf{i}}{\mathbf{i}} = 1 \\
    \text{If} \ \mathbf{i}\neq\mathbf{j}: & \quad \minor{(AB)}{\mathbf{i}}{\mathbf{j}} \geq \minor{A}{\mathbf{i}}{\mathbf{j}} + \minor{B}{\mathbf{i}}{\mathbf{j}}
\end{align*}
This proves the claim.
\end{Prf}
\end{Lem}
\begin{Rem}
  Increasing-completeness of the partial cyclic order on oriented flags is practically the same as the notion of \emph{semi-continuity from the left} in \cite{FockGoncharov}, Section 7.4. The proof we give is essentially the same, but we opted to include it since we work in the setting of oriented flags.
\end{Rem}
\begin{Prop} \label{prop:nested_tpos}
Let $g\in\PSL(n,\bR)$. Then $g$ is totally nonnegative if and only if $g\ival{F_e}{F_{w_0}} \subset \ival{F_e}{F_{w_0}}$. Moreover, $g$ is totally positive if and only if $\overline{g\sfival} \subset \sfival$.
\end{Prop}
\begin{Prf}
We first show how the inclusion statements imply total nonnegativity/positivity. By \autoref{lem:std_triple}, elements of $\ival{F_e}{F_{w_0}}$ are characterized by the fact that they admit representatives in $\SL(n,\bR)$ such that all left--bound minors (i.e. minors using the first $k$ columns) are positive. Let $F\in\ival{F_e}{F_{w_0}}$ and $M\in\SL(n,\bR)$ be such a representative for $F$. Since $gF\in\ival{F_e}{F_{w_0}}$, there is a lift $\widehat{g} \in \SL(n,\bR)$ of $g$ such that $\widehat gM$ has positive left--bound minors. By the Cauchy--Binet formula, we have
\begin{equation}    \label{eq:nested_tpos}
    \minor{(\widehat gM)}{i+1\ldots i+k}{1\ldots k} = \sum\limits_{\mathbf{j} \in \mind{k}{n}} \minor{\widehat g}{i+1\ldots i+k}{\mathbf{j}} \minor{M}{\mathbf{j}}{1\ldots k}, \qquad i+k \leq n,
\end{equation} 
so the sum on the right hand side must be positive. We will show that for any fixed $\mathbf{j}_0 \in \mind{k}{n}$, there exists a totally positive $M$ such that $\minor{M}{\mathbf{j}_0}{1\ldots k}$ is arbitrarily large compared to the other left--bound size $k$ minors. A totally positive matrix in particular has positive left--bound minors and thus represents a flag in $\ival{F_e}{F_{w_0}}$. Therefore, \eqref{eq:nested_tpos} will allow us to conclude that all minors $\minor{\widehat g}{i+1\ldots i+k}{\mathbf{j}}$ are nonnegative. This is a sufficient conditions for total nonnegativity of $\widehat g$ \cite[Proposition 2.7]{Pinkus}.

Observe that since totally positive matrices are dense in totally nonnegative matrices \cite[Theorem 2.6]{Pinkus}, it is enough to find a totally nonnegative matrix $M$ with this property. Define the submatrix $\submat{M}{\mathbf{j}_0}{1\ldots k}$ to be any totally positive $(k\times k)$--matrix and let all the other entries of $M$ be zero. Then all other left--bound size $k$ minors vanish and $M$ is totally nonnegative, so it fits our criteria. This finishes the implication ``inclusion $\Rightarrow$ totally nonnegative''.

Now assume that $\overline{g\sfival} \subset \sfival$. Then, since $gF_e,gF_{w_0} \in \sfival$, the minors
\[ \minor{\widehat g }{i+1\ldots i+k}{1\ldots k} \qquad \text{and} \qquad \minor{(\widehat g w_0)}{i+1\ldots i+k}{1\ldots k} = \minor{\widehat g }{i+1\ldots i+k}{n-k+1 \ldots n} \]
must be strictly positive. By \cite[Proposition 2.5]{Pinkus}, this is sufficient to conclude that $\widehat g$ is totally positive.

Conversely, let $g\in\PSL(n,\bR)$ be totally nonnegative and $F\in\sfival$. Let $\widehat g \in\SL(n,\bR)$ be the totally nonnegative lift and $M\in\SL(n,\bR)$ a representative for $F$. Up to replacing $M$ by $-M$, all of its left--bound minors are positive by \autoref{lem:std_triple}. Among the minors $\minor{\widehat g }{i+1\ldots i+k}{\mathbf{j}}, \ \mathbf{j}\in\mind{k}{n}$, there must be one which is strictly positive since $\widehat g$ is nonsingular. Therefore, \eqref{eq:nested_tpos} shows that left--bound connected minors of $\widehat g M$ are positive, so $gF \in \sfival$.

If $g$ is totally positive, \cite[Theorem 2.10]{Pinkus} states that it admits a decomposition as $g = LDU$, where $L \in \ltup$, $U \in \utup$ and $D$ is diagonal with positive diagonal entries. Then, $U$ stabilizes $F_e$. Moreover, $w_0 U w_0$ is lower triangular and the signs of its minors satisfy the hypotheses of \autoref{Cor:standard_opposite_triple}, so $w_0 U w_0 B_0 \in \ival{F_{w_0}}{F_e}$ which implies that $U F_{w_0} \in \ival{F_e}{F_{w_0}}$. $D$ stabilizes $\sfival$ since it fixes the endpoints. $L$ sends $F_e$ into $\sfival$ by \autoref{lem:std_triple}, fixes $F_{w_0}$ and thus sends $DU(F_{w_0})$ into $\ival{L(F_e)}{F_{w_0}}$. In particular, $(F_e,g(F_e),g(F_{w_0}),F_{w_0})$ is a cycle, so $\overline{g\ival{F_e}{F_{w_0}}} \subset \ival{F_e}{F_{w_0}}$ by properness of the partial cyclic order.
\end{Prf}
The connection to cycles observed at the end of the proof gives another characterization of totally positive matrices. Recall that $A=\stab(F_e,F_{w_0})\subset \PSL(n,\bR)$ denotes the subgroup of all diagonal matrices with positive diagonal entries.
\begin{Lem} \label{lem:cycle_tpos}
$g\in\PSL(n,\bR)$ is totally positive if and only if $(F_e,gF_e,gF_{w_0},F_{w_0})$ is a cycle. The induced map
\begin{align*} \{ g\in\PSL(n,\bR) \mid g \ &\text{totally positive} \}/A\\ &\longrightarrow \{ (F,F') \in \left(\oFlag{n}\right)^2 \mid  (F_e,F,F',F_{w_0}) \ \text{is a cycle} \}
\end{align*}
is a bijection.
\end{Lem}
\begin{Prf}
If $g$ is totally positive, we saw in the proof of the previous proposition that $(F_e,gF_e,gF_{w_0},F_{w_0})$ is a cycle. Conversely, assume that $(F_e,gF_e,gF_{w_0},F_{w_0})$ is a cycle. We will construct an element $g'\in\PSL(n,\bR)$ such that $g'(F_e) = g(F_e)$, $g'(F_{w_0}) = g(F_{w_0})$, and $g'$ is totally positive. Then $g$ and $g'$ can only differ by right--multiplication with an element of $A$, which preserves total positivity. 

Let $L\in\ltup$ be a representative for $g(F_e)$. Then $L(F_e) = g(F_e)$ and $L$ fixes $F_{w_0}$. Since $g(F_{w_0}) \in \ival{L(F_e)}{F_{w_0}}$, we have $L^{-1}g(F_{w_0}) \in \ival{F_e}{F_{w_0}}$. Let $U\in\utup$ be chosen such that $U (F_{w_0}) = L^{-1}g(F_{w_0})$. Then by the Cauchy--Binet formula, $LU$ is totally positive and we have $LU(F_e,F_{w_0}) = (g(F_e),g(F_{w_0}))$.
\end{Prf}
\begin{Cor} \label{cor:nested_tpos}
Let $I=\ival {F_2}{F_3}$ and $J=\ival {F_1}{F_4}$ be intervals in $\oFlag{n}$ such that $\overline{I} \subset J$. Then $(F_1,F_2,F_3,F_4)$ is a cycle.
\end{Cor}
\begin{Prf}
Since the action of $\PSL(n,\bR)$ on oriented transverse pairs is transitive, we may assume that $J = \sfival$. Moreover, there exists $g\in\PSL(n,\bR)$ such that $gJ = I$. By \autoref{prop:nested_tpos}, this $g$ is totally positive, so \autoref{lem:cycle_tpos} finishes the proof.
\end{Prf}

\begin{Cor} \label{cor:IntervalTotPosRepresentative}
Let $F\in \ival{F_e}{F_{w_0}}$. Then, $F=gB_0$ for some totally positive $g\in\PSL(n,\bR)$.
\end{Cor}

\begin{Lem}
  If $F$ is oriented transverse to $F'$, then there exist $F_1,F_2$ such that the quadruple $(F_1,F,F_2,F')$ is a cycle. In other words, oriented transversality characterizes those flags $F'$ in the comparable set $\comp(F)$ (\autoref{Def:comparableSet}).
  
  \begin{Prf}
    Since $F,F'$ are oriented transverse, they are in the orbit of the two standard flags by Lemma \ref{Lem:transitive_on_pairs}. By Lemma \ref{lem:std_triple} we can find $F_2\in \ival{F}{F'}$ and $F_1\in \ival{F'}{F}$, proving the claim.
  \end{Prf}
\end{Lem}
The following Lemma can be interpreted to say that if a sequence of nested intervals $\ival{F_n}{G_n}$ shrink to a single flag $F$, then the opposite intervals $\ival{G_n}{F_n}$ converge to the comparable set of $F$.
\begin{Lem}
\label{Lem:regularityOfPCO}
  Let $(F_m),(G_m) \in (\oFlag{n})^\bN$ be two sequences of oriented flags and $F\in\oFlag{n}$ satisfying the following :
  \begin{itemize}
      \item $F_m$ is increasing and converges to $F$;
      \item $G_m$ is decreasing and converges to $F$;
      \item $\cycle{F_1}{F}{G_1}$.
  \end{itemize}
  Then, 
  \[\bigcup_{m=1}^\infty \ival{G_m}{F_m} = \comp(F).\]
  \begin{Prf}
  Let $K\subset \comp(F)$ be a compact subset. We first show that there exists an $m_0$ such that $K\subset \comp(G_{m_0})\cap \comp(F_{m_0})$ for all $m\geq m_0$.
  
  Suppose not. Then there is a sequence $k_m\in K$ with $k_m \not\in \comp(F_m)$ for all $m$ or $k_m \not\in \comp(G_m)$ for all $m$; assume the first case happens. After replacing $k_m$ by a subsequence, the $i$-dimensional part $k^{(i)}_{N_m}$ is either nontransverse to $F^{(n-i)}_{N_m}$ for all $m$ or $k^{(i)}_{N_m} \oplus F^{(n-i)}_{N_m} \equalo -\bR^n$ for all $m$. Then $k^{(i)}_\infty$ must either be nontransverse to $F^{(n-i)}$ or $k^{(i)}_\infty \oplus F^{(n-i)} \equalo -\bR^n$, where $k_\infty \in K$ is any accumulation point of the sequence. This is a contradiction to $K\subset \comp(F)$.
  
  
  Now, note that by assumption $F_2$ is in the intersection $\left(\bigcup_{n=1}^\infty \ival{G_n}{F_n}\right) \cap \comp(F)$. Let $K_m$ be an exhaustion of $\comp(F)$ by connected compact sets all containing $F_2$. By the claim above, for every $m$ there is an $N$ such that $K_m\subset \comp(F_N)\cap\comp(G_N)$. But since $K_m$ is connected and $F_2\in K_m$, it must be contained in the same connected component of $\comp(G_N)\cap \comp(F_N)$ as $F_2$, which by \autoref{prop:IntervalIsConnectedComponent} is $\ival{G_N}{F_N}$. 
  \end{Prf}
\end{Lem}

\subsection{Metric on intervals}
We will define a metric on intervals of oriented flags by embedding the space of oriented flags into the product of oriented Grassmannians and then using the Pl\"ucker embedding on each factor. Let us assume that $n$ is odd for now to avoid the formal complications caused by modding out the action of $\pm 1$ in even dimension. The embeddings mentioned above take on the form

\begin{equation}    \label{eq:Plucker_embedding}
    \oFlag{n} \hookrightarrow \prod_{k=1}^{n-1} \Gro(k,n) \hookrightarrow \prod_{k=1}^{n-1} \spr \left( \bigwedge^k \bR^n \right) \cong \prod_{k=1}^{n-1} \spr \left( \bR^{\binom{n}{k}} \right),
\end{equation}

where $\Gro(k,n)$ denotes the Grassmannian of oriented $k$-planes in $\bR^n$, and $\spr$ denotes the spherical projectivization (modding out positive scalars).


Let $F_1,F_2\in\ival{F_e}{F_{w_0}}$ be a pair of oriented flags in the standard interval. This means that $F_1,F_2$ have lower triangular, unipotent, totally positive matrix representatives. In particular, all of their Pl\"ucker coordinates in the standard basis of $\Lambda^k(\bR^n)$ are positive. We conclude that the image of the above embedding lies in the product of the standard simplices in the spheres (the spherical projectivizations of the positive orthants). Since the projection from the sphere to projective space restricts to a diffeomorphism on the standard simplex, we will not distinguish between the spherical and regular projectivization from this point on. This also implies that the description applies in the same way to the case of even dimension, where we quotient by the action of $-1$. Let
\[ \iota_i \colon \oFlag{n} \to \pr\left(\bR^{\binom{n}{i}}\right) \]
denote the composition of the Pl\"ucker embedding \eqref{eq:Plucker_embedding} with the projection to the $i$-th factor.

\begin{Def}
  Define the \emph{interval distance}
  \[d_{\ival{F_e}{F_{w_0}}}(F_1,F_2):=\max_{i=1\dots n-1} d_{\Delta}(\iota_i(F_1),\iota_i(F_2)),\]
  where $d_\Delta$ is the Hilbert metric on the standard simplex. We can compute this distance explicitly using minors. If $X,Y$ are matrix representatives for $F_1,F_2$, then :
  \begin{equation*}    \label{Eq:metric_interval}
    d_{\ival{F_e}{F_{w_0}}}(F_1,F_2) = \max\limits_{\substack{1\leq k\leq n \\ \mathbf{i},\mathbf{j} \in \mathcal{I}(k,n)}} \left| \log \dbfrac{ \frac{\minor{X}{\mathbf{i}}{1\ldots k}} {\minor{X}{\mathbf{j}}{1\ldots k}} } { \frac{\minor{Y}{\mathbf{i}}{1\ldots k}} {\minor{Y}{\mathbf{j}}{1\ldots k}}  }\right| = \max_{k} \log \dbfrac{\max\limits_\mathbf{i} \left|\frac{\minor{X}{\mathbf{i}}{1\ldots k}}{\minor{Y}{\mathbf{i}}{1\ldots k}}\right| }{\min\limits_\mathbf{i} \left|\frac{\minor{X}{\mathbf{i}}{1\ldots k}}{\minor{Y}{\mathbf{i}}{1\ldots k}}\right|} .
\end{equation*} 
\end{Def}

\begin{Prop}
  The metric $d_{\ival{F_e}{F_{w_0}}}$ is invariant under the stabilizer of $\ival{F_e}{F_{w_0}}$ in $\PSL(n,\bR)$.
  \begin{Prf}
  The action of $\PSL(n,\bR)$ on $\oFlag{n}$ induces the exterior power action on each of the $\RP{\binom nk-1}$ factors. This is a projective linear action, and the stabilizer of $\ival{F_e}{F_{w_0}}$ stabilizes the standard simplex (it multiplies each Pl\"ucker coordinate by a positive number). Therefore, it acts by isometries of the Hilbert metric on each factor.
  \end{Prf}
\end{Prop}

\begin{Prop}
  \label{Prop:MetricContraction}
  Let $A\in\PSL(n,\bR)$ be such that $(F_e,AF_e,AF_{w_0},F_{w_0})$ is a cycle. Then, $A$ is a $C$-Lipschitz contraction for $d_{\ival{F_e}{F_{w_0}}}$, with $C<1$.
    \begin{Prf}
    By \autoref{lem:cycle_tpos}, a matrix $A$ such that $(F_e,AF_e,AF_{w_0},F_{w_0})$ is a cycle is totally positive. Since the $k\times k$ minors of $A$ are exactly the matrix coefficients of the exterior power $\Lambda^k A$ in the standard basis, the action on Pl\"ucker coordinates sends the positive orthant into itself. More precisely, the standard simplex is sent to another simplex, the span of the columns of $\Lambda^k A$, whose closure is contained in the standard simplex. By a classical result on Hilbert metrics (see \cite{Bir} Section 4, Lemma 1), we conclude that for each $k$, there exists $C_k<1$ such that $\Lambda^k A$ is $C_k$-Lipschitz for the Hilbert metric on the standard simplex of $\RP{\binom{n}{k-1}}$. Letting $C=\max C_k$ we get that the action of $A$ on $\ival{F_e}{F_{w_0}}$ is $C$-Lipschitz.
    \end{Prf}
\end{Prop}
\section{Schottky groups}
We now consider two types of Schottky groups which can be constructed using the partial cyclic order on oriented flags and \autoref{Def:Gen_Schottky} : Schottky groups with a purely hyperbolic model, and Schottky groups with a finite area model (see \autoref{sec:PCOSpaces}).

\subsection{Anosov representations} \label{sec:Anosov_reps}

In this section, we prove that purely hyperbolic Schottky representations in $\PSL(n,\bR)$ are Anosov with respect to the minimal parabolic $B$ if the defining intervals have disjoint closures.

First, we define a \emph{boundary map} $\xi : \partial \Gamma \rightarrow \oFlag{n}$ for any purely hyperbolic Schottky group $\rho(\Gamma)\subset \PSL(n,\bR)$. To do so, we identify $\partial \Gamma$ with the limit set of $\Gamma$ in $\RP{1}$. Then, each point $x\in \bdry$ is equal to the intersection of a unique sequence of nested intervals $I_{\gamma_k}$, where $\gamma_k$ is a reduced word of length $k$ (see \autoref{sec:PCOSpaces} for the bijection between words and intervals). Since the sequence of intervals is nested, each word $\gamma_{k+1}$ is obtained from $\gamma_k$ by adding a letter on the right, so points in $\bdry$ can also be interpreted as infinite words in the generators. Using the $k$-th order intervals $J_{\gamma_k}\subset \oFlag{n}$, we can define $\xi(x)=\bigcap_k J_{\gamma_k}$. This is well defined by the contraction of the metric on intervals from \autoref{Prop:MetricContraction}.

\begin{Prop}
\label{prop:bdrymapPurelyHyp}
The boundary map $\xi$ is continuous, equivariant, and increasing.
\begin{Prf}
Let $x_n \to x$ be a sequence in $\bdry$ converging to $x$. This implies that for any $N\in\bN$, we can find an index $n_0$ such that for all $n\geq n_0$, the first $N$ letters of $x_n$ and $x$ agree. Then $\xi(x_n)$ and $\xi(x)$ lie in the same $N$-th order interval, and by the Lipschitz contraction (\autoref{Prop:MetricContraction}) the diameter of this interval shrinks exponentially in $N$ (in the interval metric of the first order interval determined by the first letter of $x$). Thus $\xi$ is continuous.
 
 Next, we prove equivariance. Let $\gamma\in\Gamma$ be some element, expressed as a reduced word of length $l$ in the generators $A_i$ and their inverses. Then $\gamma x \in \bdry$, as an infinite sequence in the generators, is simply the concatenation of the finite word $\gamma$ and the infinite word $x$. Denoting by $w^{(j)}$ the word $w$ truncated after the first $j$ letters, the corresponding nested sequence of intervals is $I_{(\gamma x)^{(j)}}, \ j\in\bN$, and $\xi$ maps $\gamma x$ to the unique point in $\bigcap\limits_j J_{(\gamma x)^{(j)}}$. By definition, we have $\rho(\gamma)J_{x^{(j)}} = J_{\gamma x^{(j)}}$, so this intersection point agrees with $\rho(\gamma)\xi(x)$.
 
 Finally, we show that $\xi$ is increasing. Let $x,y,z \in \bdry \cong \Lambda_\Gamma$ be three points such that $\cycle xyz$. Then there are indices $K,L,M \in \bN$ such that the intervals $I_{x^{(K)}}, I_{y^{(L)}}, I_{z^{(M)}}$ are in increasing configuration. Their image intervals $J_{x^{(K)}}, J_{y^{(L)}}, J_{z^{(M)}}$ satisfy the same cyclic relations and contain the points $\xi(x),\xi(y),\xi(z)$ respectively, so we conclude $\cycle {\xi(x)}{\xi(y)}{\xi(z)}$.
\end{Prf}
\end{Prop}
To prove that Schottky representations are Anosov, we will use a characterization that was given in \cite[Theorem 1.7]{KapovichLeebPortiMorseActions}, specialized to the case at hand.
\begin{Def}[{\cite[Definition 4.1]{KapovichLeebPortiCharacterizations}}]
A sequence $(g_n) \in G^\bN$ is \emph{B--contracting} if there exist flags $F_+,F_- \in \Flag{n}$ such that
\[ g_n|_{C(F_-)} \xrightarrow{n\to\infty} F_+ \]
locally uniformly, where $C(F_-)$ denotes the set of flags transverse to $F_-$.\\
If $\Gamma$ is a group and $\rho\colon\Gamma\to G$ a homomorphism, we also call a sequence $(\gamma_n) \in \Gamma^\bN$ B--contracting if $(\rho(\gamma_n))$ is.
\end{Def}
\begin{Def}[{\cite[Definition 4.25]{KapovichLeebPortiCharacterizations}}]
Let $H<G$ be a subgroup. Then the \emph{$\Flag{n}$--limit set} consists of all flags $F_+$ as in the previous definition for all contracting sequences $(g_n) \in H^\bN$.
\end{Def}
The following is part of \cite[Theorem 1.1]{KapovichLeebPortiCharacterizations}.
\begin{Thm} \label{thm:Anosov_asymptotically_embedded}
  Let $\Gamma$ be a word hyperbolic group and $\rho\colon\Gamma\to G$ a representation. Then $\rho$ is $B$--Anosov if and only if:
  \begin{enumerate}
    \item There is a $\rho$--equivariant embedding
    \[ \xi \colon \bdry \to \Flag{n} \]
    whose image is the $\Flag{n}$--limit set of $\Gamma$ such that for any $x\neq y \in \bdry$, $\xi(x)$ and $\xi(y)$ are transverse.
    \item Every diverging sequence $\gamma_n \to \infty$ in $\Gamma$ has a $B$--contracting subsequence.
  \end{enumerate}
\end{Thm}
\begin{Prop}
  Let $\rho: \Gamma \to \PSL(n,\bR)$ be a purely hyperbolic generalized Schottky representation. Then $\rho$ is $B$-Anosov.
  \begin{Prf}
    As $\Gamma$ is a free group, it is word hyperbolic. We have constructed a boundary map into the space of oriented flags $\oFlag{n}$ (\autoref{prop:bdrymapPurelyHyp}), which we can project to flags to get a boundary map $\pi \circ \xi$ which is continuous, equivariant, and transverse.
    
    We now show how to find $B$--contracting subsequences. Let $\gamma_n$ be a divergent sequence in $\Gamma$. As $\Gamma$ is free, each $\gamma_n$ can be written in a unique way as a reduced word in the generators $A_i$ and their inverses. Let $\gamma_n^{(k)}$ denote the subword consisting of the first $k$ letters, and $\gamma_n^{(-k)}$ the subword consisting of the last $k$ letters. After replacing $\gamma_n$ by a subsequence, we may assume that $\ell(\gamma_n) > 2n$ (where $\ell$ denotes word length), and the subwords $\gamma_m^{(n)}$ and $\gamma_m^{(-n)}$ are constant for $m\geq n$; denote them by $a_n$ and $b_n$. Then $\rho(\gamma_m)$ maps the interval $-J_{b_n^{-1}}$ into the interval $J_{a_n}$ for all $m\geq n$. Note that $a_{n+1}$ resp. $b_{n+1}^{-1}$ are obtained from $a_n$ resp. $b_n^{-1}$ by adding a letter on the right. Thus $J_{b_n^{-1}}$ is a nested sequence of intervals converging to a point $\xi(x)$ and $J_{a_n}$ is a nested sequence of intervals converging to a point $\xi(y)$. By \autoref{Lem:regularityOfPCO}, the opposites $-J_{b_n^{-1}}$ form an increasing sequence converging to the comparable set $\comp(\xi(x))$. By \autoref{Cor:UniqueTransverseLift}, the transverse set $C(\pi \circ \xi(x))$ is the projection $\pi\circ \comp(\xi(x))$.
    
    This argument also shows that the image of $\pi \circ \xi$ is exactly the $\Flag{n}$--limit set of $\Gamma$.
    
    
    
  \end{Prf}
\end{Prop}
\subsection{Positive representations}   \label{sec:positive_reps}
Let $\Sigma$ be a closed, oriented surface with finitely many punctures and $\Gamma = \pi_1(\Sigma)$.

Fixing a finite area hyperbolization of $\Sigma$, the set of lifts of punctures of $\Sigma$ is called the \emph{Farey set} and denoted by $\mathscr{F}_\infty \subset \RP{1}$.

\begin{Def}
  A map $\beta:\mathscr{F}_\infty \rightarrow \Flag{n}$ is \emph{positive} if for every finite cyclically ordered tuple $p_1,\dots,p_N$ of $\mathscr{F}_\infty\subset \RP{1}$, the image tuple $\beta(p_1)\dots \beta(p_N)$ is positive (\autoref{def:positivity}).
\end{Def}

\begin{Def}[\cite{FockGoncharov}]
  A representation $\rho: \pi_1(\Sigma) \rightarrow \PSL(n,\bR)$ is \emph{positive} if it admits an equivariant, positive map $\beta : \mathscr{F}_\infty \rightarrow \Flag{n}$.
\end{Def}

\begin{Lem}
  The positive map $\beta$ associated to a positive representation $\rho$ can be lifted to an increasing, equivariant map $\hat{\beta} : \mathscr{F}_\infty \rightarrow \oFlag{n}$. 
  \begin{Prf}
    Fix $p_0\in \mathscr{F}_\infty$ and a lift $\hat\beta(p_0)\in \oFlag{n}$ of $\beta(p_0)$. By \autoref{Cor:UniqueTransverseLift}, for every $p\neq p_0$ in $\mathscr{F}_\infty$, there exists a unique lift of $\beta(p)$ to oriented flags which is oriented transverse to $\hat\beta(p_0)$, and we define $\hat\beta(p)$ to be this lift. It remains to show that $\hat\beta$ is increasing and equivariant.
    
    Let $x,y,z\in \mathscr{F}_\infty$ such that $(p_0,x,y,z)$ is a cycle. By positivity of $\beta$, this means that $\beta(p_0),\beta(x),\beta(y),\beta(z)$ can be represented as
    \[B,gB,ghB,w_0B\]
    where $g,h\in\ltup$. Without loss of generality, the lift we chose for $\beta(p_0)$ is $B_0$. Since $w_0 B_0$, $gB_0$, and $ghB_0$ are oriented transverse to $B_0$, they are the cosets representing $\hat\beta(x),\hat\beta(y),\hat\beta(z)$. They form a cycle in $\oFlag{n}$, and so in particular $\cycle{\hat\beta(x)}{\hat\beta(y)}{\hat\beta(z)}$ showing that $\hat\beta$ is increasing. This also shows that $\hat\beta(q)$ is oriented transverse to $\hat\beta(p)$ whenever $q\neq p$, which will use in the proof of equivariance below. 
    
    To show equivariance, let $\gamma\in\Gamma$ and $p_\gamma \in \RP{1}$ a fixed point of $\gamma$. By \cite{FockGoncharov}, the eigenvalues of $\rho(\gamma)$ are all positive. Moreover, Fock and Goncharov show that we can extend $\beta$ to the fixed points of elements in $\Gamma$ preserving positivity (Theorem 7.2). If we extend $\beta$ this way, $\rho(\gamma)$ preserves every lift of $\beta(p_\gamma)$ by positivity of all eigenvalues, and so in particular it preserves the unique lift $\hat\beta(p_\gamma)$ which is oriented transverse to $\hat\beta(p_0)$. For $q\in\mathscr{F}_\infty$ different from $p_\gamma$, since $\PSL(n,\bR)$ preserves oriented transversality we have that $\rho(\gamma)\hat\beta(q)$ is oriented transverse to $\hat\beta(p_\gamma)$. Moreover, it is a lift of $\rho(\gamma)\beta(q)=\beta(\gamma q)$, and so it must be equal to $\hat\beta(\gamma q)$. 
  \end{Prf}
\end{Lem}

\begin{Thm}
  A representation $\rho : \Gamma \rightarrow \PSL(n,\bR)$ is positive if and only if it admits a Schottky presentation.
\begin{Prf}
Let $\rho : \Gamma \rightarrow \PSL(n,\bR)$ be a Schottky representation with a finite area model (\autoref{Def:Gen_Schottky}). By \autoref{Thm:SchottkyLeftContinuousMap}, it admits an equivariant left-continuous and increasing boundary map $\xi:\RP{1}\rightarrow \oFlag{n}$. Denote by $\pi :\oFlag{n} \rightarrow \Flag{n}$ the projection from oriented flags to flags. Then, the map $(\pi \circ \xi)|_{\mathscr{F}_\infty}$ is a positive equivariant map and so $\rho$ is a positive representation.

Conversely, if $\rho$ is positive, then it admits an equivariant positive map $\beta :\mathscr{F}_\infty \rightarrow \Flag{n}$, which we can lift to an equivariant increasing map $\hat{\beta} : \mathscr{F}_\infty \rightarrow \oFlag{n}$. The finite area hyperbolization of $\Sigma$ used in defining $\mathscr{F}_\infty$ can be given a Schottky presentation where the intervals have endpoints in $\mathscr{F}_\infty$. To do this, cut the surface along geodesic arcs beginning and ending at a puncture until the complement of the arcs is a topological disk. The edge identifications of this disk will give a Schottky presentation as claimed. We can push this presentation to the space of oriented flags using the map $\hat{\beta}$, associating an interval $(a,b)$ to the interval $\ival{\hat{\beta}(a)}{\hat{\beta}(b)}$. By equivariance and increasingness, this is a Schottky presentation of $\rho$.
\end{Prf}
\end{Thm}
\section{Fundamental Domains}
\label{sec:fundamentalD}
In this section, we use the techniques developed above in order to build fundamental domains bounded by piecewise linear faces for the action of positive representations on the sphere and projective space. The space and the domains vary with the value of $n \mod 4$. In particular, there is no cocompact domain of discontinuity in either the sphere or projective space when $n=4k+1$ (this will be treated in \cite{Stecker}).

\subsection{Fundamental domains in projective space}
Let $n=2k$ for $k\in \mathbb{N}$. We will construct a fundamental domain for $\PSL(n,\bR)$ Schottky groups in projective space $\RP{n-1}$.

To do so, we associate to each oriented transverse pair of complete oriented flags $F_1,F_2\in\oFlag{n}$ a \emph{halfspace} in $\RP{n-1}$. These halfspaces will play the same role as half planes in the hyperbolic plane when constructing classical Schottky groups.

\begin{Def}
  The \emph{open halfspace} associated to an oriented tranverse pair $F_1,F_2$ is :
  \[\Half(F_1,F_2) := \{[v]\in \RP{n-1} ~|~ \Var^+_\mathcal{E}(v)\leq k-1\}\]
  Similarly, the \emph{closed halfspace} of $F_1,F_2$ is:
  \[\overline{\Half}(F_1,F_2) := \{[v]\in \RP{n-1} ~|~ \Var^-_\mathcal{E}(v)\leq k-1\},\]
  where in both cases we take the variation with respect to a basis $\mathcal{E}$ adapted to $F_1,F_2$ (see \autoref{def:variation}).
\end{Def}

The action of $\PSL(n,\bR)$ on halfspaces is given by
\[ g\Half(F_1,F_2) = \Half(gF_1,gF_2). \]
The following lemma is easily checked:

\begin{Lem}
  $\Half(F_1,F_2)$ is open. The closure of $\Half(F_1,F_2)$ is $\overline{\Half}(F_1,F_2)$.
\end{Lem}

\begin{figure}
    \centering
    \includegraphics[width=.8\textwidth]{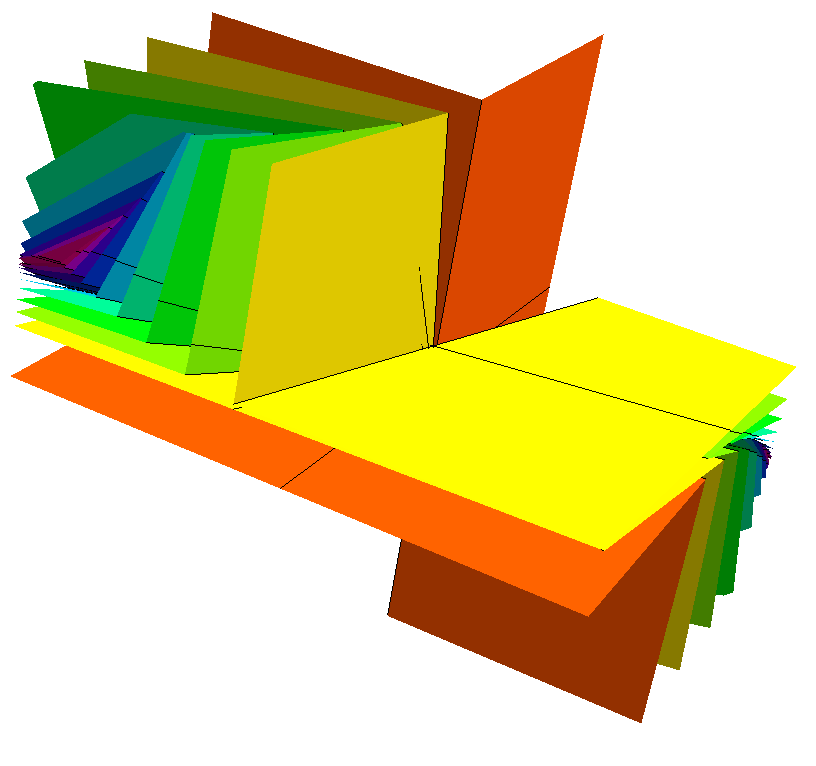}
    \caption{A family of nested halfspaces in $\RP{3}$.}
    \label{fig:NestedHalfspacesRP3}
\end{figure}

\begin{Ex}
  For $n=4$, the standard halfspace is the closure of four orthants defined by sign vectors $(1,1,1,1),(-1,1,1,1),(-1,-1,1,1),(-1,-1,-1,1)$. \autoref{fig:NestedHalfspacesRP3} depicts a family of halfspaces in $\RP{3}$. The largest halfspace in the nested family is the standard halfspace.
\end{Ex}

\begin{Lem} \label{Lem:HalfSpaceComplementEven}
The complement of the open halfspace $\Half(F_1,F_2)$ is the closed opposite halfspace, $\overline{\Half}(F_2,F_1)$.
\begin{Prf}
  If $\mathcal{E}=(e_1,\dots,e_n)$ is a basis adapted to $F_1,F_2$, then
  \[\mathcal{E}'=(e_n,-e_{n-1},\dots,e_2,-e_1)\]
  is adapted to $F_2,F_1$. The action of $w_0$ on coordinates changes coordinates from $\mathcal{E}$ to $\mathcal{E}'$.
  
  Upper and lower variation are related in the following way : $\Var^+_\mathcal{E}(v) = n-1-\Var_\mathcal{E}^-(w_0 v)$, which means that $\Var^+_\mathcal{E}(v)\geq k$ if and only if $\Var^-_\mathcal{E}(w_0 v)\leq k-1$.
\end{Prf}
\end{Lem}

\begin{Lem} \label{NestedEvenHalfspaces}
  Let $A$ be an $n\times n$ totally positive matrix. Then,
  \[A (\overline{\Half}(F_e,F_{w_0})) \subset \Half(F_e, F_{w_0}).\]
  \begin{Prf}
    The Lemma follows immediately from the variation diminishing property of totally positive matrices given in \autoref{Thm:VariationDiminishing}.
  \end{Prf}
\end{Lem}

We now have the prerequisites to prove the main disjointness theorem for halfspaces.

\begin{Thm}\label{Thm:HalfspaceDisjointness}
  Let $F_1,F_2,F_3,F_4$ be a cyclically ordered quadruple of oriented flags in $\bR^n$ with $n$ even. Then, $\Half(F_1,F_2)$ and $\Half(F_3,F_4)$ are disjoint.
  \begin{Prf}
    We normalize so that $F_2 = F_e$ and $F_1 = F_{w_0}$. Since $F_e,F_3,F_4,F_{w_0}$ are in cyclic order, by \autoref{lem:cycle_tpos} there exists a totally positive matrix $A$ with $AF_e=F_3$ and $AF_{w_0}=F_4$. Then, by \autoref{NestedEvenHalfspaces}, $\Half(F_3,F_4)\subset \Half(F_e,F_{w_0})$, and by \autoref{Lem:HalfSpaceComplementEven} $\Half(F_e,F_{w_0})\subset \Half(F_{w_0},F_e)^c$ which finishes the proof. 
  \end{Prf}
\end{Thm}

\begin{Rem}
  \autoref{Thm:HalfspaceDisjointness} is not an equivalence. It is possible to construct disjoint halfspaces which are not defined by cyclically ordered oriented flags.
\end{Rem}

The following technical lemma will be useful in proving an equivalent definition of a halfspace which is more closely related to intervals in full flags. Recall that $\ltup \subset \PSL(n,\bR)$ denotes the subgroup of lower triangular, unipotent, totally positive matrices.

\begin{Lem} \label{Lem:VariationSpanTotallyPositive}
  If $\Var^+(v) < m$, there exists $M\in\ltup$ such that $v$ is in the span of the first $m$ columns of $M$. Moreover, if $\Var^+(v)=m-1$, we can write the linear combination such that the coefficient of the $m$-th column of $M$ has the same sign as that used for the last coordinate of $v$ when determining $S^+(v)$.
  \begin{Prf}
  Let $u=Jv$ where $J$ is a diagonal matrix with entries alternating betwen $1$ and $-1$. Then, $\Var^-(u)=n-1-\Var^+(v)>n-m-1$. We will show that there is a matrix $M\in\ltup$ such that the last $n-m$ entries of $Mu$ vanish. This implies that $u$ is in the span of the first $m$ columns of $M^{-1}$, and so $v$ is in the span of the first $m$ columns of $JM^{-1}J$, which is triangular totally positive (as in the proof of \autoref{Lem:involution}).
  
  Denote by $x_i(a)$ the lower triangular matrix with $1$s on the diagonal and only one other nonzero entry $a$ in position $(i+1,i)$. We also denote
  \[C_r(a_1,\dots,a_{r}) = x_{r}(a_{r})x_{r-1}(a_{r-1})\dots x_1(a_1).\]
  The action of $x_i(a)$ on a vector is to add a multiple of the $i$-th coordinate to the $(i+1)$-th coordinate.
  
  By the decomposition theorem (\autoref{Thm:TotPosDecomposition}), a product of the form
  \[C_1 \dots C_{n-2}C_{n-1}\]
  is in $\ltup$ if $C_i = C_i(a_1,\dots,a_{i})$ with $a_j>0$.
  
  We claim that we can choose strictly positive values for $a_1,\dots,a_{n-1}$ so that the last entry of $C_{n-1} u$ vanishes, and $\Var^-(C_{n-1} u) = \Var^-(u)-1$.
    
     To prove this, assume that the last sign change in $u$ is between $u_i$ and $u_{i+1}$. One can choose small enough positive values for $a_1,\dots, a_{i-1}$ so that the signs of $u_1,\dots,u_{i}$ do not change under the action of $x_{i-1}(a_{i-1})\dots x_1(a_1)$ (if one or more of these values is zero, then it will get the sign of the previous entry, which does not change $\Var^-$). By assumption, $u_{i+1}\dots u_n$ all have the same sign, but different from $u_i$, so we can choose $a_i,\dots a_{n-2}>0$ to make $(C_{n-1}(a_1,\dots,a_{n-2})u)_{i+1}, \dots, (C_{n-1}(a_1,\dots,a_{n-2})u)_{n-1}$ have the same sign as $u_i$, and we can choose $a_{n-1}$ in order to make $(C_{n-1} u)_n$ vanish, unless $u_n$ is zero. If $u_n=0$, then use the same strategy to make the last nonzero entry of $u$ map to $0$ under $C_{n-1}$, and so all the following entries will stay zero. Since the only sign changes which could matter for $\Var^-$ were in the entries after $i$, where the sign change between $u_i$ and $u_{i+1}$ was removed, $\Var^-(C_{n-1} u)=\Var^-(u)-1$. By induction, we can make the last $m$ entries vanish with a product of the form
    \[C_1 C_2\dots C_{n-2}C_{n-1}.\]
    For the last part of the statement, we now assume that $\Var^+(v)=m-1$. Note that by the argument above we found $M\in\ltup$ such that the last $n-m$ entries of $Mu$ vanish, but also we have that $\Var^-(Mu)=0$ since we lowered the variation of $u$ by one every time we multiplied by one of the $C_r$. Therefore, $\Var^+(JMu) = n-1$ and the first $m$ entries of $JMu$ are nonzero and alternate in sign, which means $\Var^-(JMu)=m-1$. We have
    \[v=(JM^{-1}J)(JMu),\]
    where $JM^{-1}J\in \ltup$. Let $U\in \utup$, and write
    \[v=(JM^{-1}JU)(U^{-1}JMu) = P u'.\]
    Note that the matrix $P=JM^{-1}JU$ is (strictly) totally positive. By the variation diminishing theorem for totally nonnegative matrices \cite[Theorem 3.4]{Pinkus},
    \[\Var^-(u') = \Var^-(U^{-1}JMu)\geq \Var^-(UU^{-1}JMu) = m-1.\]
    By upper triangularity of $U^{-1}$, the last $n-m$ coordinates of $u'$ are also zero and the $m$th coordinate of $u'$ is equal to that of $JMu$. This forces $\Var^-(u')=m-1$.
    
    Since $\Var^+(v) =\Var^+(Pu')= \Var^-(u')=m-1$, we can use the equality case of the variation diminishing theorem (\autoref{Thm:VariationDiminishing}). It says precisely that the sign of the last coordinate of $v$ used in determining $S^+(v)$ must be the same as that of the last nonzero coordinate of $u'$, which is also the same as the last coordinate of $JMu$.
  \end{Prf}
\end{Lem}

\begin{Prop}\label{Prop:HalfspacesAndIntervals}
  Halfspaces and intervals are related in the following way:
  \[\Half(F_1,F_2) = \bigcup_{F\in \ival{F_1}{F_2}} \mathbb{P}F^{(k)}\]
  \begin{Prf}
    It is sufficient to prove this for $(F_1,F_2)=(F_e,F_{w_0})$. Then, the statement of the proposition is equivalent to:
    
    A vector $v\in \bR^n$ satisfies $\Var^+(v)\leq k-1$ if and only if it is in the span of the first $k$ columns of some totally positive, lower triangular, unipotent $n\times n$ matrix $M$.
    
    We first show that the right hand side is included in the halfspace. Let $P$ be a totally positive representative for some flag $F\in\ival{F_e}{F_{w_0}}$ (it exists by \autoref{cor:IntervalTotPosRepresentative}). Then, $v=P^{(k)}u$ where $P^{(k)}$ is the matrix formed by the first $k$ columns of $P$, and $u\in \bR^k$. By \autoref{Thm:VariationDiminishing},
    \[S^+(v) = S^+(P^{(k)}u)\leq S^-(u)\leq k-1.\]
    
    For the other inclusion, we want to show that any $v\in\bR^n$ with $\Var^+(v)<k$ is in the span of the first $k$ columns of a matrix $M\in\ltup$. This follows from \autoref{Lem:VariationSpanTotallyPositive} with $m=k$.
  \end{Prf}
\end{Prop}


With this notion of halfspace, we can construct fundamental domains as follows:
\begin{Thm}
\label{thm:FundamentalD}
Let $\Gamma$ be a Schottky group in $\PSL(2k,\bR)$ defined by the disjoint, cyclically ordered intervals $\ival{F_1}{G_1}, \dots, \ival{F_{2g}}{G_{2g}}$. Then, the set
  \[D = \RP{2k-1}-\bigcup_{i=1}^{2g} \Half(F_i,G_i) = \bigcap_{i=1}^{2g} \overline{\Half}(G_i,F_i)\]
  is a fundamental domain for the action of $\Gamma$ on its maximal domain of discontinuity in projective space.
  \begin{Prf}
    Let $\gamma\in\Gamma$ be a nontrivial element. Then, $\gamma D$ and $D$ have disjoint interiors by the cyclic ordering of the Schottky intervals and \autoref{Thm:HalfspaceDisjointness}.
    
    We now show that
    \begin{equation}    \label{eq:fundamental_domain_orbit}
        \bigcup_{\ell(\gamma)\leq m} \gamma D = \RP{2k-1} - \bigcup_{\ell(\gamma)=m+1} \Half(J_\gamma)  ,
    \end{equation} 
    where $J_\gamma$ is the interval associated to $\gamma$ and $\Half(J_\gamma)$ is the halfspace associated to this interval (see \autoref{sec:PCOSpaces} for the bijection between words and intervals).
    
    Let $\gamma$ be a nontrivial element. We can write $\gamma=\gamma' a$, where $a$ is the last letter of $\gamma$ (so $\ell(\gamma')=\ell(\gamma)-1$). For every generator or inverse of a generator such that $a'\neq a^{-1}$, 
    \[\gamma \Half(J_{a'}) = \Half(J_{\gamma a'}).\]
    Moreover,
    \[\gamma \Half(J_{a^{-1}}) = \Half(\gamma' a J_{a^{-1}}) = \Half(\gamma' (-J_a)) = \Half(-J_{\gamma}).\]
    We conclude that 
    \[ \gamma \bigcup_{i=1}^{2g} \Half(F_i,G_i) = \Half(-J_\gamma) \cup \bigcup_{\gamma'\in L_{+1}(\gamma)} \Half(J_{\gamma'}), \]
    where $L_{+1}(\gamma) \subset \Gamma$ contains the $2g-1$ elements of length $\ell(\gamma) + 1$ of the form $\gamma a'$. Thus
    \[ \gamma D = \RP{2k-1} - \gamma \bigcup_{i=1}^{2g} \Half(F_i,G_i) = \overline{\Half}(J_\gamma) - \bigcup_{\gamma'\in L_{+1}(\gamma)} \Half(J_{\gamma'}). \]
    Moreover, since $\ell(\gamma)$-th order intervals are cyclically ordered,
    \[ \overline{\Half}(J_\gamma) \subset \RP{2k-1} - \bigcup_{\stackrel{\ell(\gamma')=\ell(\gamma)}{\gamma' \neq \gamma}} \Half(J_{\gamma'}) \subset \RP{2k-1} - \bigcup_{\stackrel{\ell(\gamma')=\ell(\gamma)+1}{\gamma' \not\in L_{+1}(\gamma) }} \Half(J_{\gamma'}), \]
    so the inclusion ``$\subset$'' in \eqref{eq:fundamental_domain_orbit} follows.\\
    Conversely, let $x\in \RP{2k-1} - \bigcup_{\ell(\gamma)=m+1} \Half(I_\gamma)$ and let $m(x)\geq 0$ be minimal such that $x$ is contained in $\RP{2k-1} - \bigcup_{\ell(\gamma)=m(x)+1} \Half(I_\gamma)$. If $m(x)=0, x\in D$. Now assume that $m(x)>0$. Then 
    \begin{align*}
        x & \in \left( \RP{2k-1} - \bigcup_{\ell(\gamma)=m(x)+1} \Half(I_\gamma) \right) - \left( \RP{2k-1} - \bigcup_{\ell(\gamma)=m(x)} \Half(I_\gamma) \right) \\
        & = \bigcup_{\ell(\gamma) = m(x)} \overline{\Half}(I_{\gamma}) - \bigcup_{\ell(\gamma)=m(x)+1} \Half(I_\gamma),
    \end{align*} 
    so there exists $\gamma\in\Gamma$ such that $\ell(\gamma)=m(x)$ and $x\in\overline{\Half}(I_\gamma) - \bigcup_{\gamma'\in L_{+1}(\gamma)} \Half(I_{\gamma'}) = \gamma D$. 
    
    Now by Lipschitz contraction of intervals (\autoref{Prop:MetricContraction}) and \autoref{Prop:HalfspacesAndIntervals}, 
    \[ \bigcap_{m=1}^\infty \bigcup_{\ell(\gamma)=m} \Half(I_\gamma) = \bigcup_{x\in\bdry} \xi(x)^{(k)}, \]
    so the full domain $\Gamma D$ is the open and dense set
    \[ \RP{2k-1} - \bigcup_{x\in\bdry} \xi(x)^{(k)}. \]
  \end{Prf}
\end{Thm}

\begin{figure}
    \centering
    \includegraphics[width=.8\textwidth]{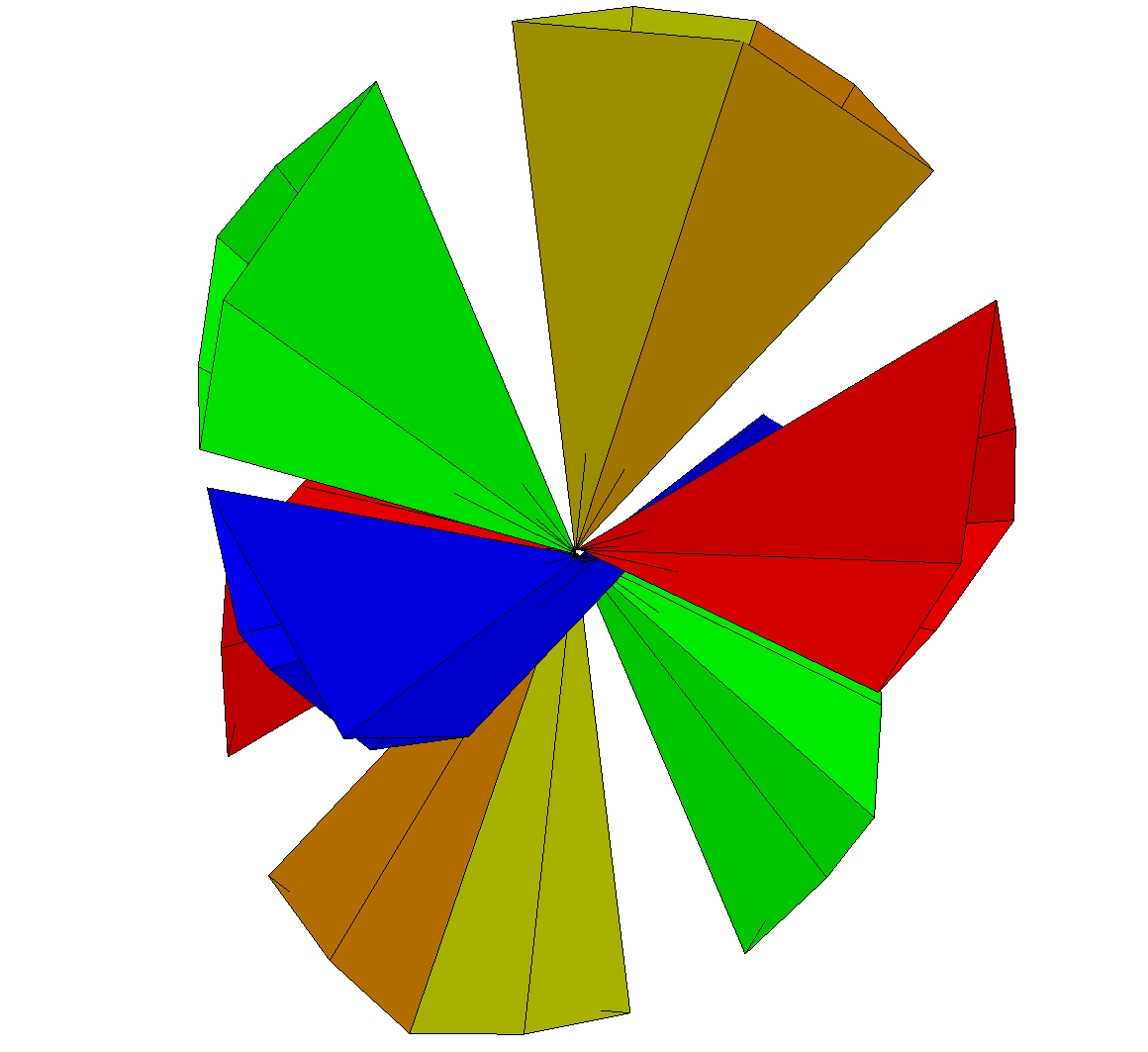}
    \caption{Four disjoint $\RP{3}$ halfspaces bounding a fundamental domain.}
    \label{fig:FourDisjointRP3}
\end{figure}

\subsection{Fundamental domains in the sphere}
For the rest of this section, let $n=4k+3$. We will construct a fundamental domain in the sphere $S^{4k+2}$. The construction is analogous to the even dimensional case, but one has to be slightly more careful with the notion of halfspace. The reason for this is that a vector in $\bR^{4k+3}$ can have anywhere from $0$ to $4k+2$ sign changes, and so we cannot define a halfspace to be ``all vectors with less that half the maximum number of sign changes''. The case $\Var^+(v)=2k+1$ has to be ``split in half''. 

\begin{Def}
  Let $F_1,F_2\in\oFlag{4k+3}$ be a pair of oriented transverse flags. Let $\mathcal{E}$ be a basis adapted to this pair. The \emph{open halfspace} in the sphere $S^{4k+2}$ associated to $(F_1,F_2)$ is the set
 \begin{align*}
     \Half(F_1,F_2) = \{[v]\in S^{4k+2} :& \Var^+_\mathcal{E}(v)\leq 2k+1 \text{ and in case of equality,}\\ &\text{the sign of the last coordinate of $v$}\\&\text{used in determining }\Var^+_\mathcal{E}(v)\text{ is positive}\}
  \end{align*}
  The \emph{closed halfspace}, similarly, is the set
  \begin{align*}
     \overline{\Half}(F_1,F_2) = \{[v]\in S^{4k+2} :& \Var^-_\mathcal{E}(v) \leq 2k+1 \text{ and in case of equality,}\\ &\text{the last nonzero coordinate of $v$ is positive}\}
  \end{align*}
\end{Def}

In this definition we could have chosen to ask that, in the case of equality, the last nonzero coordinate is negative instead. Indeed, there are two possible conventions for halfspaces and they give fundamental domains for the two different possible maximal domains of discontinuity in the sphere.

\begin{figure}
    \centering
    \includegraphics[width=.5\textwidth]{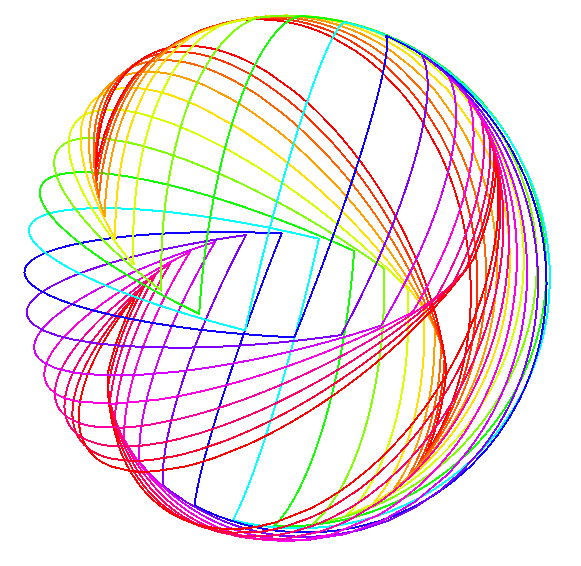}
    \caption{A 1-parameter family of nested halfspaces in $S^2$ (only the boundaries are shown).}
    \label{fig:NestedHalfspaces}
\end{figure}

\begin{Ex}
  For $n=3$, the standard halfspace is the closure of four orthants defined by sign vectors $(1,1,1),(-1,1,1),(-1,-1,1),(-1,-1,-1)$.
\end{Ex}
As in the even dimensional case, totally positive matrices preserve the standard halfspace. \autoref{fig:NestedHalfspaces} depicts a 1-parameter family of totally positive matrices acting on a halfspace.

\begin{Lem}
  Let $A$ be an $n \times n$ totally positive matrix. Then,
  \[A(\overline{\Half}(F_e,F_{w_0}))\subset \Half(F_e,F_{w_0})\]
  
  \begin{Prf}
    This is a direct consequence of the variation diminishing property in \autoref{Thm:VariationDiminishing}.
  \end{Prf}
\end{Lem}

\begin{figure}
    \centering
    \includegraphics[width=.5\textwidth]{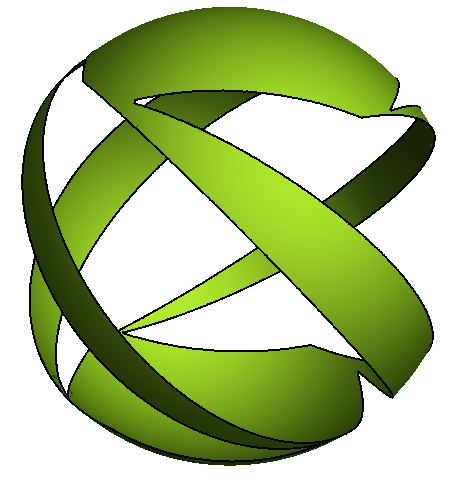}
    \caption{A fundamental domain in $S^2$ defined by the intersersection of four halfspaces.}
    \label{fig:FundamentalDomainSphere}
\end{figure}

The analog of \autoref{Prop:HalfspacesAndIntervals} in this setting is that a halfspace in the sphere is a union of \emph{positive hemispheres}.
\begin{Def}
Let $F\in\oFlag{4k+3}$. The $k$-dimensional positive hemisphere associated to $F$ is the subset of the sphere
\[F^{(k)}_+ = \{[v]\in S^{4k+2} ~|~ F^{(k-1)} \oplus [v] \equalo F^{(k)}\}.\]
\end{Def}
\begin{Prop}
Halfspaces and intervals in spheres are related by
\[\Half(F_1,F_2) = \bigcup_{F\in\ival{F_1}{F_2}} F^{(2k+2)}_+.\]
\begin{Prf}
The proof is analogous to that of \autoref{Prop:HalfspacesAndIntervals}. For the standard pair $F_e,F_{w_0}$, the statement translates to the equivalence between the two following conditions:

\begin{itemize}
    \item $S^+(v)\leq 2k+1$ and if $S^+(v)=2k+1$ then the sign of the last coordinate of $v$  used in determining $S^+(v)$ is positive,
    \item $v$ is in the span of the first $2k+2$ columns of some $L\in\ltup$ and the coefficient of the $(2k+2)$nd column is positive.
\end{itemize}

If $v$ satisfies the second condition, we have $v=Lu$ for $u\in\bR^{4k+3}$, $u_{2k+2}>0$ and $u_j=0$ for $j>2k+2$. Let $U\in\utup$ so that $LU$ is totally positive and $LUB_0=LB_0$. Since $U$ is upper triangular, the $2k+2$-th coordinate of $U^{-1}u$ is that same as that of $u$ and all coordinates after vanish. Using the variation diminishing theorem (\autoref{Thm:VariationDiminishing}),
\[S^+(v) = S^+((LU)(U^{-1}u)) \leq S^-(U^{-1}u) \leq 2k+1.\]
If $S^+(v)=2k+1$, this is an equality, and the variation diminishing theorem then tells us that the sign of the last coordinate of $v$ used in determining $S^+(v)$ is the same as the sign of the last nonzero coordinate of $U^{-1}u$, which is positive.

The converse follows directly from \autoref{Lem:VariationSpanTotallyPositive} with $m=2k+2$.
\end{Prf}
\end{Prop}
\begin{Thm}
  Let $\Gamma$ be a Schottky group in $\PSL(4k+3,\bR)$ defined by the disjoint, cyclically ordered intervals $\ival{F_1}{G_1}, \dots, \ival{F_{2g}}{G_{2g}}$. Then, the set
  \[D = S^{4k+2}-\bigcup_{i=1}^{2g} \Half(F_i,G_i) = \bigcap_{i=1}^{2g} \overline{\Half}(G_i,F_i)\]
  is a fundamental domain for the action of $\Gamma$ on its domain of discontinuity in the sphere.
  
  \begin{Prf}
    The proof is completely analogous to that of \autoref{thm:FundamentalD}.
  \end{Prf}
\end{Thm}
\appendix
\section{Anti-de Sitter crooked planes}
In this appendix we will show that our notion of halfspace in $\RP{3}$ coincides with that of \emph{anti-de Sitter crooked halfspace} introduced in \cite{DGKArcComplex} and studied in \cite{DGKAntideSitter}, \cite{GoldmanCrookedSurf}. More precisely, an $\AdS^3$ crooked halfspace is the restriction to the projective model of $\AdS^3$ of an $\RP{3}$ halfspace as defined in \autoref{sec:fundamentalD}.

\begin{Def}
The $3$-dimensional anti-de Sitter space $\AdS^3$ is the group $\PSL(2,\bR)$ endowed with the Lorentzian metric given by its Killing form.
\end{Def}

The isometry group of $\AdS^3$ is $\PSL(2,\bR)\times \PSL(2,\bR)$ acting by left and right multiplication.

\begin{Def}
Given a geodesic $\ell\subset \bH^2$ in the hyperbolic plane, the associated $\AdS^3$ \emph{right crooked plane} $C(\ell)$ is the set of isometries $g\in\PSL(2,\bR)$ which have a nonattracting fixed point on $\overline{\ell}\subset \overline{\bH^2}$.
\end{Def}

We will consider the following embedding of $\AdS^3$ in $\RP{3}$:
\[ \iota : \AdS^3 \rightarrow \RP{3} \]
\[ \begin{pmatrix} a & b\\c & d\\ \end{pmatrix} \mapsto \begin{bmatrix}\frac{a+d}{2}\\-b\\ \frac{a-d}{2}\\ c\end{bmatrix}.\]

Its image is the projectivization of the set of negative vectors for the signature $(2,2)$ quadratic form $q(v) = -v_1^2 + v_3^2 - v_2 v_4$.

\begin{Prop}
The standard halfspace $\Half$ in $\RP{3}$ intersects the image of $\iota$ in one of the two crooked $\AdS^3$ halfspaces bounded by $C(\ell_0)$, where $\ell_0$ is the geodesic represented by the positive imaginary axis in the Poincar\'e upper half plane.

\begin{Prf}
The boundary of the standard halfspace is the closure of the projectivization of vectors which have signs $(+,0,+,*),(*,+,0,+),(0,+,*,-)$, or $(+,*,-,0)$.

Let us analyse each of these cases. If $M=\begin{pmatrix} a & b\\c & d\end{pmatrix}$ and $\iota(M)$ has signs $(+,0,+,*)$, then $b=0$ and so $M$ fixes $0$. Moreover, $a>0$ and $-a<d<a$ which means that $0$ is a repelling fixed point.

Similarly, if $\iota(M)$ has signs $(+,*,-,0)$, then $c=0$ and $M$ fixes $\infty$. Now, $d>0$ and $-d<a<d$ so $\infty$ is repelling.

If $\iota(M)$ has signs $(*,+,0,+)$, then $a=d$, $b<0$ and $c>0$. This means that $M$ is an elliptic element fixing the point $i\sqrt{\frac{-b}{c}}\in i \bR$.

Finally, the set with signs $(0,+,*,-)$ does not intersect the image of $\iota$ since any vector with these signs is positive for the quadratic form $q$.

Conversely, any $M\in\PSL(2,\bR)$ in the crooked plane must fall into one of these categories or their closure (if $M$ is parabolic or the identity).
\end{Prf}
\end{Prop}

\bibliography{bibliography}
\bibliographystyle{alpha}
\end{document}